\documentclass[11pt]{article}
\usepackage{graphicx,amsmath,amsfonts,amssymb,bbm,amsthm}
\usepackage{psfrag}
\usepackage{epsf}
\usepackage{hyperref}
\makeatletter\@addtoreset {equation}{section}\makeatother

\theoremstyle{plain}
\newtheorem{theo}{Theorem}[section]

\newtheorem{cor}[theo]{Corollary}
\theoremstyle{remark}
\newtheorem{rem}[theo]{Remark}

{\hspace*{\fill}$\rule{.3\baselineskip}{.35\baselineskip}$\end{trivlist}}

\renewcommand{\geq}{\geqslant}
\renewcommand{\leq}{\leqslant}

\renewcommand{\phi}{\varphi}
\newcommand{\be}{\begin{eqnarray}}
\newcommand{\ee}{\end{eqnarray}}

\newcommand{\eps}{\varepsilon}

\begin{document}

\title{\bf On the orbital stability of Gaussian solitary waves \\ in the log-KdV equation}

\author{R\'emi Carles$^{1}$ and Dmitry Pelinovsky$^{2}$ \\
{\small $^{1}$ CNRS \& Univ. Montpellier, Math\'ematiques (I3M),
34095 Montpellier, France} \\
{\small $^{2}$ Department of Mathematics, McMaster
University, Hamilton, Ontario, Canada, L8S 4K1}  }

\date{\today}
\maketitle

\begin{abstract}
We consider the logarithmic Korteweg--de Vries (log--KdV) equation, which models solitary waves
in anharmonic chains with Hertzian interaction forces. By using
an approximating sequence of global solutions of the regularized generalized
KdV equation in $H^1(\mathbb{R})$ with conserved $L^2$ norm and energy,
we construct a weak global solution of the log--KdV equation in
a subset of $H^1(\mathbb{R})$.
This construction yields conditional orbital stability of Gaussian solitary waves of the log--KdV equation,
provided uniqueness and continuous dependence of the constructed solution holds.

Furthermore, we study the linearized log--KdV equation at the Gaussian solitary wave
and prove that the associated linearized operator has a purely discrete spectrum
consisting of simple purely imaginary eigenvalues
in addition to the double zero eigenvalue. The eigenfunctions, however, do
not decay like Gaussian functions but have algebraic decay.
Using numerical approximations, we show that
the Gaussian initial data do not spread out
but produce visible radiation at the left slope of the Gaussian-like pulse in the time
evolution of the linearized log--KdV equation.
\end{abstract}

\section{Introduction}

Solitary waves in anharmonic chains with Hertzian interaction forces are modelled by the Fermi--Pasta--Ulam (FPU)
lattices with non-smooth nonlinear potentials \cite{Nesterenko}. Recently, the FPU lattice equations
in the limit of small anharmonicity of the Hertzian interaction forces
were reduced to the following logarithmic Korteweg--de Vries (log-KdV) equation \cite{Chat,JP13}:
\begin{equation}
\label{logKdV}
v_t + v_{xxx} + (v \log|v|)_x = 0, \quad (x,t) \in \mathbb{R} \times \mathbb{R}.
\end{equation}
Here and in what follows,  the subscripts denote the partial derivatives.

The log--KdV equation (\ref{logKdV}) has a two-parameter family of
Gaussian solitary waves
\begin{equation}
\label{soliton-orbit}
v(x,t) = e^{c} v_G(x-ct-a), \quad a,c \in \mathbb{R},
\end{equation}
where $v_G$ is a symmetric standing wave given by
\begin{equation}
\label{Gaussian}
v_G(x) := e^{\frac{1}{2} - \frac{x^2}{4}}, \quad x \in \mathbb{R}.
\end{equation}

The ultimate goal of this work is to prove the nonlinear orbital stability of Gaussian
solitary waves (\ref{soliton-orbit}) in the log--KdV equation (\ref{logKdV}). The main problem is, of course,
the limited smoothness of the log--KdV equation, where the nonlinearity
$f(v) = v \log|v|$ is continuous but not differentiable at $v = 0$,
whereas the energy
\begin{equation}
\label{energy}
E(v) = \frac{1}{2} \int_{\mathbb{R}} \left[ (v_x)^2 - v^2 \log|v| \right] dx + \frac{1}{4} \int_{\mathbb{R}} v^2 dx
\end{equation}
is not a $C^2$ functional at $v = 0$.

Although $E(v)$ is not $C^2$ at $v = 0$, the second variation of $E(v)$ at $v_G$ is well
determined by the Schr\"{o}dinger operator with a harmonic potential
\begin{equation}
\label{Schrodinger-operator}
L := - \partial_x^2 + \frac{1}{4} (x^2 - 6).
\end{equation}
Note that the spectrum of $L$ in $L^2(\mathbb{R})$ consists of
equally spaced simple eigenvalues
$$
\sigma(L) = \{ -1,0,1,2,\ldots\},
$$
which include exactly one negative eigenvalue with the eigenvector
$v_G$ (defined without normalization). Therefore, $E(v)$ is not convex at $v_G$ in a subspace of
$H^1(\mathbb{R})$. Nevertheless, the second variation of $E(v)$ at $v_G$
given by $E_c(u) = \frac{1}{2} \langle L u, u \rangle_{L^2}$ is positive in the constrained space
\begin{equation}
\label{constrained-space}
X_c := \left\{ u \in H^1(\mathbb{R})\cap L^2_1(\mathbb{R}) : \quad
\langle v_G, u \rangle_{L^2} = 0 \right\},
\end{equation}
where $L^2_1(\mathbb{R}) := \{ u \in L^2(\mathbb{R}) : \; x u \in L^2(\mathbb{R})\}$.
At the linearized approximation, this constraint fixes $\| v \|_{L^2}^2$ at $\| v_G \|_{L^2}^2$.
Based on these facts, the linear orbital stability
of the Gaussian solitary wave $v_G$ can be deduced in the following sense.

Consider the time evolution of the linearized log-KdV equation at the Gaussian wave $v_G$,
\begin{equation}
\label{linlogKdV}
\left\{ \begin{array}{l} u_t = \partial_x L u, \quad t > 0, \\
u |_{t = 0} = u_0. \end{array} \right.
\end{equation}
Note that the quadratic energy function $E_c(u)$ is constant in time for
smooth solutions of the linearized log--KdV equation (\ref{linlogKdV}).
We say that the Gaussian solitary wave $v_G$ is linearly orbitally stable in $H^1(\mathbb{R})$,
if for every $u_0 \in X_c$, there exists a unique global solution $u$
of the linearized log--KdV equation (\ref{linlogKdV}) in a subspace of $L^{\infty}(\mathbb{R},H^1(\mathbb{R}))$
which satisfies the following bound
\begin{equation}
\label{global-bound}
\| u(\cdot,t) \|_{H^1} \leq C(\| u_0 \|_{H^1},\| u_0 \|_{L^2_1}), \quad t \in \mathbb{R},
\end{equation}
for some $t$-independent positive constant $C$ that depends on the initial norms
$\| u_0 \|_{H^1}$ and $\| u_0 \|_{L^2_1}$.

The following theorem is based on the fact that the conserved
quantity $E_c(u)$ is positive if $u \in X_c$
and controls the squared $H^1$ norm of the solution $u$ in time $t$.
Although this theorem was not formulated in \cite{JP13}, it can be deduced
from the arguments developed in this work, which rely on symplectic projections
and the energy method. Due to the scaling invariance (\ref{soliton-orbit}),
the same stability result holds for all values of parameter $a$ and $c$.

\begin{theo}\cite{JP13}
Gaussian solitary wave $v_G$ of the log--KdV equation (\ref{logKdV}) is linearly
orbitally stable in $H^1(\mathbb{R})$. \label{theorem-GP}
\end{theo}

To extend this result to the proof of nonlinear orbital stability of Gaussian
solitary waves, we need local and global well-posedness of the log--KdV equation
(\ref{logKdV}). In a similar context of the log--NLS equation, global well-posedness
and orbital stability of the Gaussian solitary waves were proved
by Cazenave and Haraux \cite{Caz1} and Cazenave \cite{Caz2} (these
results were later summarized in the monograph \cite[Section 9.3]{Caz}).

The idea from \cite{Caz1,Caz2} is to approximate the logarithmic nonlinearity
by a smooth nonlinearity, to construct a sequence of global bounded solutions
of the regularized system in $C(\mathbb{R},H^1(\mathbb{R}))$, and then to prove convergence
of a subsequence to the limit, which solves the original log--NLS
equation in a weak sense.
Uniqueness of the weak solution does not come in this formalism for free, but it
can be proved with an additional trick involving log-nonlinearity
\cite[Lemma 9.3.5]{Caz},
which is based on the bound (\ref{log-bound}) below. Similar ideas
were recently used by Carles and Gallo \cite{CG} for analysis of the
NLS equation with power-like
nonlinearity and a sublinear damping term (see also \cite{CaOz-p} for
various extensions).

While trying to adopt the programme above to the log--KdV equation (\ref{logKdV}), we come
to two main difficulties. The first one is that local solutions of the generalized
KdV equation with smooth nonlinearities exist in the space $C([-t_0,t_0],H^s(\mathbb{R}))$
for $s > \frac{3}{2}$ (see \cite{Kato,Zhidkov} for two independent proofs of these results).
Therefore, to push $s$ to lower values, in particular, to $s = 1$,
we need to adopt the formalism of Kenig, Ponce, and Vega \cite{KPV93}, which was originally developed to
the KdV equation with integer powers. The second difficulty is that the proof
of uniqueness requires us to consider more restrictive solutions of the log--KdV equation
than the ones constructed in the proof of existence. As a result, the solutions
we establish in a subspace of $H^1(\mathbb{R})$ have non-decreasing
$L^2$ norm and energy. Unlike the case of weak
solutions of the NLS equation (e.g., treated in \cite{GV}),
we cannot establish even the $L^2$ conservation in the log--KdV
equation (\ref{logKdV}). See Remark \ref{remark-L2} below for further details.

To make sense of the energy (\ref{energy}) of the log--KdV equation
(\ref{logKdV}), we shall work with functions in the class
\begin{equation}
\label{energy-space}
X := \left\{ v \in H^1(\mathbb{R}) : \quad v^2 \log|v| \in L^1(\mathbb{R}) \right\}.
\end{equation}
The following theorem gives the main result on the existence of weak solutions
of the log--KdV equation (\ref{logKdV}) in the energy space $X$.

\begin{theo}
For any $v_0 \in X$, there exists a global solution
$v \in L^{\infty}(\mathbb{R},X)$ of the log--KdV equation (\ref{logKdV})
such that
\begin{equation}
\| v(t) \|_{L^2} \leq \| v_0 \|_{L^2}, \quad E(v(t)) \leq E(v_0), \quad \mbox{\rm for all} \;\; t \in \mathbb{R}.
\end{equation}
Moreover, if $\partial_x \log|v| \in L^{\infty}((-t_0,t_0) \times \mathbb{R})$,
then the solution $v$ exists in $C((-t_0,t_0),X)$,
is unique for every $t \in (-t_0,t_0)$,
depends continuously on the initial data $v_0 \in X$, and satisfies
conservation of $\| v(t) \|_{L^2}$ and
$E(v(t))$ for all $t \in (-t_0,t_0)$.
\label{theorem-CP}
\end{theo}

Unfortunately, $\partial_x \log|v|$ is unbounded as $|x| \to \infty$
for the Gaussian solitary wave $v_G$. Therefore, it is unclear
from Theorem \ref{theorem-CP} if uniqueness and continuous dependence hold for such Gaussian solutions.
Convexity of the energy functional $E(v)$ at $v_G$ and
global well-posedness of the Cauchy problem for
the log--KdV equation (\ref{logKdV}) in $X$ are the two main ingredients of the
nonlinear orbital stability of the Gaussian solitary wave $v_G$ in $H^1(\mathbb{R})$. Because of limitations
in Theorem \ref{theorem-CP}, we can only obtain the conditional nonlinear orbital stability,
where the condition is that the global solution $v \in L^{\infty}(\mathbb{R},X)$
of the log--KdV equation (\ref{logKdV}) constructed in
Theorem \ref{theorem-CP} is unique and depends continuously on the initial data $v_0 \in X$.

To define the nonlinear orbital stability, we use the standard theory (see \cite{Pava} for review of
this theory). We say that the Gaussian solitary wave $v_G$ is orbitally stable in $H^1(\mathbb{R})$
if for any $\eps > 0$ there exists $\delta > 0$ such that
for any $v_0 \in X \subset H^1(\mathbb{R})$ satisfying
\begin{equation}
\| v_0 - v_G \|_{H^1} \leq \delta,
\end{equation}
there exists a unique solution $v \in C(\mathbb{R},X)$
of the log--KdV equation (\ref{logKdV}) satisfying
\begin{equation}
\inf_{a \in \mathbb{R}} \| v(\cdot,t) - v_G(\cdot + a) \|_{H^1} \leq \eps,
\end{equation}
for all $t \in \mathbb{R}$. The following conditional orbital stability
result follows from Theorems \ref{theorem-GP} and \ref{theorem-CP} with
the standard arguments \cite{Pava}.

\begin{cor}
Gaussian solitary wave $v_G$ of the log--KdV equation (\ref{logKdV}) is
orbitally stable in $H^1(\mathbb{R})$ under condition that
the solution $v \in L^{\infty}(\mathbb{R},X)$ in Theorem \ref{theorem-CP}
is unique and depends continuously on the initial data $v_0 \in X$.
\label{theorem-orbital}
\end{cor}

Since Corollary \ref{theorem-orbital} does not give a proper nonlinear orbital stability result,
we shall also look at the orbital stability problem from a different point of view.
We first inspect properties of the linearized operator $\partial_x L$
in the linearized log--KdV equation (\ref{linlogKdV}). The following theorem
gives the spectral stability of the Gaussian solitary wave $v_G$
with precise characterization of eigenvalues and eigenvectors of the
linear operator $\partial_x L$ in $L^2(\mathbb{R})$.

\begin{theo}
The spectrum of $\partial_x L$ in $L^2(\mathbb{R})$ is purely discrete and
consists of a double zero eigenvalue and a symmetric sequence of simple
purely imaginary eigenvalues $\{ \pm i \omega_n \}_{n \in \mathbb{N}}$
such that $0 < \omega_1 < \omega_2 < ...$ and $\omega_n \to \infty$ as $n \to \infty$.
The double zero eigenvalue corresponds to the Jordan block
\begin{equation}
\label{null-space}
\partial_x L \partial_x v_G = 0, \quad \partial_x L v_G = - \partial_x v_G,
\end{equation}
whereas the purely imaginary eigenvalues $\lambda = \pm i \omega_n$ correspond to the
eigenfunctions $u = u_{\pm n}(x)$, which are smooth in $x$ but decay
algebraically as $|x| \to \infty$.
\label{theorem-spectral}
\end{theo}

\begin{rem}
Because the spectrum of $\partial_x L$ is purely discrete by Theorem \ref{theorem-spectral},
we have no chance to assume dispersive decay estimates near the Gaussian solitary wave $v_G$.
This indicates that no asymptotic stability result can hold for solitary waves in
the log--KdV equation (\ref{logKdV}). \label{remark-1}
\end{rem}

\begin{rem}
Remark \ref{remark-1} agrees with the result of \cite[Proposition 4.3]{Caz2} stating
that the $L^p$ norms at the solution $v$ for any $p \geq 2$ including $p = \infty$
may not vanish as $t \to \infty$ (or in a finite time), hence the solution does not scatter
to zero. Although this statement was proved for the log--NLS equation, it is based on the consideration
of the same energy functional $E(v)$ as in the log--KdV equation (\ref{logKdV}).
This non-scattering result is related
to the fact that the $L^p$ norm at the family (\ref{soliton-orbit})
can be scaled to be arbitrarily small by using the scaling parameter $c \in \mathbb{R}$.
\label{remark-2}
\end{rem}

Theorem \ref{theorem-spectral} can be used to provide an alternative proof of Theorem \ref{theorem-GP}.
On the other hand, because of the algebraic decay of the eigenfunctions in Theorem \ref{theorem-spectral},
any function of $x$ that decays like the Gaussian function as $|x| \to \infty$ cannot be simply represented as series
of eigenfunctions of the linearized operator $\partial_x L$. To explain the importance of such
representations, we set $v(x,t) := v_G(x) + w(x,t)$ and obtain the equivalent log--KdV equation
\begin{equation}
\label{logKdV-w}
w_t = \partial_x L w - \partial_x N(w),
\end{equation}
where the nonlinear term $N(w)$ is given by
$$
N(w) := w \log\left( 1 + \frac{w}{v_G} \right) + v_G \left[ \log \left( 1 + \frac{w}{v_G} \right) - \frac{w}{v_G} \right].
$$
It is clear that the nonlinear term $N(w)$ does not behave uniformly in $x$ unless
$w$ decays at least as fast as $v_G$ in (\ref{Gaussian}). On the other hand,
if $w(x,t) = v_G(x) h(x,t)$, where $h$ is a bounded function in its variables, then
$N(w) = v_G n(h)$, where $n(h) := h \log(1 + h) + \log(1 + h) - h$ is analytic in $h$ for any $h \in (-1,1)$.
This observation on the nonlinear term $N(w)$ inspires us to consider
solutions of the linearized log--KdV equation (\ref{linlogKdV}) starting with
the initial data $u_0 \in X_c$ such that $u_0 = v_G h_0$ for $h_0 \in
L^{\infty}(\mathbb{R}) \cap L^2(\mathbb{R})$.

We have undertaken numerical simulations of this linear Cauchy problem
to illustrate that solutions $u(x,t)$ of the linearized log--KdV equation (\ref{linlogKdV})
with Gaussian initial data do not spread out as the time variable evolves.
Nevertheless, they produce visible radiation at the left slope of the Gaussian solutions.
Further studies are needed to figure out if the nonlinear
orbital stability of the Gaussian solitary wave $v_G$ can be proved in the framework of
the nonlinear evolution problem (\ref{logKdV-w})  with initial data $w_0 = v_G h_0$
for small $h_0 \in L^{\infty}(\mathbb{R}) \cap L^2(\mathbb{R})$.

The rest of the paper consists of the following. Theorem~\ref{theorem-CP}
is proved in Section~2. Theorem~\ref{theorem-spectral} is proved in Section 3.
Section 4 contains the proof of Theorem~\ref{theorem-GP} (which is alternative to the one given in \cite{JP13})
and the numerical simulations of the time evolution
of the linearized log--KdV equation (\ref{linlogKdV}) starting with Gaussian initial data.

\section{Global solutions of the log--KdV equation}

To prove Theorem~\ref{theorem-CP}, we follow the algorithm developed in
\cite{Caz1,Caz2} for the log--NLS equation (see also \cite[Section 9.3]{Caz}),
which shares common features with the approach of Ginibre and Velo
\cite{GV}. We slightly simplify the functional framework from
\cite{Caz1,Caz2},
by using Strichartz estimates instead of arguments coming from convex
analysis. In the context of KdV equations, these tools
were developed by Kenig, Ponce, and Vega \cite{KPV93}. 

The strategy of the proof is the following. We first
regularize the logarithmic nonlinearity near the origin and construct a sequence of
approximating solutions via a contraction principle as in \cite{KPV93}.
Then, we derive uniform estimates on the $H^1(\mathbb{R})$
norm of the approximating solutions. This allows us to pass to the limit and
obtain a global solution $v$ in the energy space $X$ with non-increasing $\| v \|_{L^2}$ and energy $E(v)$.
In the last subsection, we establish uniqueness of the global solutions $v$
of the log--KdV equation (\ref{logKdV}) with bounded $\partial_x \log|v|$
by using special properties of the logarithmic nonlinearity.

\subsection{Approximating solutions}

We shall regularize the behavior of the logarithmic nonlinearity $f(v) = v \log|v|$
near $v = 0$. Compared to the approximation considered in \cite{Caz1,Caz2}, we
will work with a different (much simpler) approximation of the logarithmic nonlinearity.

For any fixed $\eps > 0$, let us define the family of regularized  nonlinearities in the form
\begin{equation}
\label{approximation}
f_{\eps}(v) = \left\{ \begin{array}{l} f(v), \quad \; |v| \geq \eps, \\
p_{\eps}(v), \quad |v| \leq \eps, \end{array} \right.
\end{equation}
where $p_{\eps}$ is an odd polynomial of degree $2m+1$ such that
$\partial^k_v p_{\eps}(\eps) = \partial_v^k f(\eps)$ for all $0 \leq k \leq m$.
In this way, we construct $f_{\eps} \in C^m(\mathbb{R})$ for any $m \in \mathbb{N}$.
For instance, for $m = 1$, we can compute
\begin{equation}
\label{polynomial}
p_{\eps}(v) = \left( \log(\eps) - \frac{1}{2} \right) v + \frac{1}{2 \eps^2} v^3,
\end{equation}
which yields $f_{\eps}(v) \in C^1(\mathbb{R})$. For $m = 2$, we compute
\begin{equation}
\label{polynomial-2}
p_{\eps}(v) = \left( \log(\eps) - \frac{3}{4} \right) v + \frac{1}{\eps^2} v^3 - \frac{1}{4 \eps^4} v^5,
\end{equation}
which yields $f_{\eps}(v) \in C^2(\mathbb{R})$.  As pointed out in Remark \ref{remark-space} below,
for the proof of Theorem~\ref{theorem-CP}, it is actually sufficient to consider the
approximation $f_{\eps}(v)$ with $m=2$.

Global behavior of the function $f_{\eps}(v)$ is still determined by the logarithmic nonlinearity
$f(v)$. If $m \geq 1$, the function $f_{\eps}(v)$ is globally Lipschitz and for any fixed $\eps > 0$
there exists a positive constant $C_{\eps}$ such that
\begin{equation}
\label{bound-1}
|f_{\eps}(v) - f_{\eps}(u)| \leq C_{\eps} (|v| + |u|) |v-u|, \quad {\rm for \;\; every} \;\; v,u \in \mathbb{R}.
\end{equation}
If $m \geq 2$, we also have
\begin{equation}
\label{bound-2}
|f'_{\eps}(v) v_x - f'_{\eps}(u) u_x | \leq C_{\eps} \left( (|v|+|u|) |v_x - u_x|
+ (|v_x|+|u_x|) |v-u| \right),
\end{equation}
for every $v,u,v_x,u_x \in \mathbb{R}$, where $C_{\eps}$ is another constant,
which may change from one line to another line. Of course,
we realize from examples (\ref{polynomial}) and (\ref{polynomial-2})
that $C_{\eps} \to \infty$ as $\eps \to 0$.

For a given initial data $v_0$, we shall now consider a sequence
of the approximating Cauchy problems associated with the
generalized KdV equations
\begin{equation}
\label{generalizedKdV}
\left\{ \begin{array}{l} v^{\eps}_t + v^{\eps}_{xxx} + f'_{\eps}(v^{\eps}) v^{\eps}_x  = 0, \quad t > 0, \\
v^{\eps} |_{t = 0} = v_0, \end{array} \right.
\end{equation}
Using the linear estimates and the contraction principle
from the work of Kenig, Ponce and Vega \cite{KPV93}, we have the following
local well-posedness result.

\begin{theo}
\label{theorem-KPV93}
Fix $s > \frac{3}{4}$ and assume that $f_{\eps} \in C^2(\mathbb{R})$
satisfy the global Lipschitz estimates (\ref{bound-1}) and (\ref{bound-2}).
For any $v_0 \in H^s(\mathbb{R})$, there exists a time $T(\| v_0 \|_{H^s}) > 0$ and a unique
solution of the generalized KdV equation (\ref{generalizedKdV}) satisfying
\begin{eqnarray*}
& (1) & v^{\eps} \in C([-T,T],H^s(\mathbb{R})), \\
& (2) & v^{\eps}_x \in L^4([-T,T],L^{\infty}(\mathbb{R})), \\
& (3) & \| D_x^s v^{\eps}_x \|_{L^{\infty}_x L^2_T} < \infty, \\
& (4) & \| v^{\eps} \|_{L^2_x L^{\infty}_T} < \infty,
\end{eqnarray*}
where
$$
\| f \|_{L^p_x L^q_T} := \left( \int_{\mathbb{R}} \left( \int_{-T}^T |f(t,x)|^q dt \right)^{p/q} dx \right)^{1/p}.
$$
Moreover, the solution $v^{\eps}$ depends continuously on the initial data $v_0$ in $H^s(\mathbb{R})$ and
\begin{equation}
\label{L2-conservation}
\| v(t) \|_{L^2} = \| v_0 \|_{L^2} ,\quad \mbox{\rm for \; every} \;
t \in [-T,T].
\end{equation}
If in addition $m\geq 2$ and $s\geq 1$, then the energy is conserved:
\begin{equation}
E_{\eps}(v^{\eps}(t)) = E_{\eps}(v_0),\quad \mbox{\rm for \; every} \;
t \in [-T,T],
\label{energy-conservation}
\end{equation}
where
\begin{equation}
\label{energy-regularized}
E_{\eps}(v) := \frac{1}{2} \int_{\mathbb{R}}  (v_x)^2 dx -
\int_{\mathbb{R}} W_{\eps}(v) dx,\quad W_{\eps}(v) := \int_0^v f_{\eps}(v) dv.
\end{equation}
\end{theo}

\begin{rem}
Although Theorem \ref{theorem-KPV93} was proved for the KdV equation with $f(v) = v^2$ in \cite{KPV93},
Duhamel's principle was used to write the Cauchy problem (\ref{generalizedKdV}) in the integral form
\begin{equation}
\label{integral-KdV}
v^{\eps}(0) = S(t) v_0 - \int_0^t S(t-t') \left( f'_{\eps}(v^{\eps}(t')) v^{\eps}_x(t') \right) dt',
\end{equation}
where $S(t)$ is the solution operator associated with the group $e^{-t \partial_x^3}$. Contraction principle
is now applied in the same spirit as in \cite{KPV93} provided the regularized nonlinearity $f_{\eps}$
satisfies the global Lipschitz estimates (\ref{bound-1}) and (\ref{bound-2}).
\end{rem}

\begin{rem}
The solution is extended globally for any integer $s \in \mathbb{N}$ thanks
to the a priori energy estimates provided that the nonlinearity $f_{\eps}$ is
sufficiently smooth. We only use this construction for $s = 1$, when
$f_{\eps} \in C^2(\mathbb{R})$ is sufficient. In this context, the approximation
(\ref{approximation}) with (\ref{polynomial-2}) can be used for the rest of this work.
\label{remark-space}
\end{rem}

\subsection{Uniform energy estimates}

We work with local solutions of Theorem \ref{theorem-KPV93} for $s = 1$. Besides
conservation of the $L^2$ norm in (\ref{L2-conservation}), the generalized KdV equation (\ref{generalizedKdV})
admits conservation of the energy in (\ref{energy-conservation}), as it is stated in Theorem \ref{theorem-KPV93}.

For the purpose of construction of global solutions in $H^1(\mathbb{R})$,
we need to bound $E_{\eps}(v)$ from below by the squared $H^1(\mathbb{R})$ norm.
To do so, we only need to bound positive values of $W_{\eps}(v)$ from above.

Since $f'(v) = \log|v| + 1 \to -\infty$ as $v \to 0$,
we realize that $f_{\eps}'(v) = p'_{\eps}(v) \ll -1$ for $|v| \leq \eps$ 
if $\eps > 0$ is sufficiently small. Indeed, it follows from
the explicit example (\ref{polynomial-2})
that 
$$
f'_{\eps}(v) = p_{\eps}'(v) = \log(\eps) + \mathcal{O}(1) \quad \mbox{\rm for all} \;\; |v| \leq \eps 
\quad \mbox{\rm as} \;\; \eps \to 0.
$$
Therefore, if $\eps > 0$ is sufficiently small, we have
\begin{equation}
\label{relation-W-small}
W_{\eps}(v) = \frac{1}{2} \left[ \log(\eps) + \mathcal{O}(1) \right] v^2 \leq 0, \quad |v| \leq \eps.
\end{equation}
Let us define
\begin{equation}
\label{potential-W}
W(v) := \frac{1}{2} v^2 \log|v| - \frac{1}{4} v^2,
\end{equation}
and denote the positive part of $W$ by $[W]_+$.
Note that $W(v) \geq 0$ for $|v| \geq \sqrt{e}$.
Because $f_{\eps}(v) = f(v)$ for $|v| \geq \eps$,
there exists a positive constant $C_m$, which only depends on $m$ such that
\begin{equation}
\label{relation-W}
W_{\eps}(v) = W(v) + C_m \eps^2, \quad |v| \geq \eps.
\end{equation}
For instance, integration of the explicit example (\ref{polynomial-2}) yields  $C_2 = \frac{1}{12}$.
It follows from (\ref{relation-W-small}), (\ref{potential-W}), and (\ref{relation-W}) that there exists
a positive $\eps$-independent constant $C$  such that
\begin{equation}
\label{bound-potential}
[W_{\eps}(v)]_+ \leq [W(v)]_+ + C_m v^2 \leq C |v|^3, \quad {\rm for \;\; every} \;\; v \in \mathbb{R}.
\end{equation}
By Sobolev's embedding of $H^1(\mathbb{R})$ into $L^{\infty}(\mathbb{R})$,
we obtain from (\ref{energy-regularized}) and (\ref{bound-potential})
\begin{eqnarray*}
E_{\eps}(v) & \geq & \frac{1}{2} \| v_x \|^2_{L^2} - C \| v \|_{L^{\infty}} \| v \|_{L^2}^2 \\
& \geq & \frac{1}{2} \| v_x \|^2_{L^2} - C \| v \|_{H^1} P(v),
\end{eqnarray*}
where $P(v) := \frac{1}{2} \| v \|_{L^2}^2$. Because $E_{\eps}(v^{\eps})$ and $P(v^{\epsilon})$
are constants of motion for the solution $v^{\epsilon}$ of the generalized KdV
equation (\ref{generalizedKdV}) in Theorem \ref{theorem-KPV93} with $s = 1$, we obtain the time-uniform
bound on the $H^1$ norm of $v^{\eps}$:
\begin{equation}
\label{bound-solution}
\| v^{\eps}(t) \|_{H^1} \leq C P(v_0) + \sqrt{2 E_{\eps}(v_0) + 2 P(v_0) + C^2 P^2(v_0)} < \infty, \quad t \in [-T,T].
\end{equation}
By using bound (\ref{bound-solution}), we use a standard continuation argument for solutions of the integral equation
(\ref{integral-KdV}) obtained via a contraction principle and continue the local solution $v^{\eps} \in C([-T,T],H^1(\mathbb{R}))$ to a global
solution $v^{\eps} \in C(\mathbb{R},H^1(\mathbb{R}))$. The global solution satisfies the time-uniform bound
(\ref{bound-solution}) extended for every $t \in \mathbb{R}$.

Additionally, because of (\ref{bound-potential}), there exists an $\eps$-independent
constant $C$ such that
\begin{eqnarray}
\nonumber
\| [W_{\eps}(v^{\eps}(t))]_+ \|_{L^1} & \leq & C \|v^{\eps}(t)\|_{L^{\infty}} \| v^{\eps}(t) \|_{L^2}^2 \\
& \leq & C P(v_0) \left[ C P(v_0) + \sqrt{2 E_{\eps}(v_0) + 2 P(v_0) + C^2 P^2(v_0)}\right].
\label{bound-solution-2}
\end{eqnarray}
For the negative part of $W_{\eps}$, we obtain
\begin{eqnarray}
\nonumber
\| [W_{\eps}(v^{\eps}(t))]_- \|_{L^1} & = & \| [W_{\eps}(v^{\eps}(t))]_+ \|_{L^1} - \int_{\mathbb{R}} W_{\eps}(v^{\eps}(t)) dx \\
\nonumber
& = & \| [W_{\eps}(v^{\eps}(t))]_+ \|_{L^1} + E_{\eps}(v^{\eps}(t)) - \frac{1}{2} \| v^{\eps}_x \|_{L^2}^2 \\
& \leq & C P(v_0) \left[ C P(v_0) + \sqrt{2 E_{\eps}(v_0) + 2 P(v_0) + C^2 P^2(v_0)}\right] + E_{\eps}(v_0).
\label{bound-solution-3}
\end{eqnarray}
Bounds (\ref{bound-solution-2}) and (\ref{bound-solution-3}) allow us
to control $\| W_{\eps}(v^{\eps}(t))\|_{L^1}$ for every $t \in \mathbb{R}$.

\subsection{Passage to the limit}
\label{sec:pass-lim}

We shall now consider the limit $\eps \to 0$ for the sequence of global approximating solutions
$v^{\eps} \in C(\mathbb{R},H^1(\mathbb{R}))$ satisfying the generalized KdV equations (\ref{generalizedKdV}).
We recall definition (\ref{energy-space})
of the energy space $X$ for the log--KdV equation (\ref{logKdV}).
Assume that $v_0 \in X$ so that $E(v_0) < \infty$. Clearly $E_{\eps}(v_0) < \infty$ for any $v_0 \in H^1(\mathbb{R})$.

Since $f_{\eps}(0) = f(0) = 0$, we have the pointwise limit $f_{\eps}(v) \to f(v)$
as $\eps \to 0$ for every $v \in \mathbb{R}$
for the sequence of regularized nonlinearities in (\ref{approximation}). Consequently, we have
the pointwise limit of $W_{\eps}(v)$ to the potential $W(v)$ given by (\ref{potential-W}):
$$
W_{\eps}(v) \to W(v) \quad \mbox{\rm as} \quad \eps \to 0 \quad \mbox{\rm for \; every} \;\; v \in \mathbb{R}.
$$
In view of  (\ref{polynomial-2}), there exists a positive $\eps$-independent constant $C$ such that
\begin{equation*}
\left| W_\eps(v)- \frac{1}{2} (\log \eps) v^2\right|\leq Cv^2, \quad
  |v|\leq \eps.
\end{equation*}
We infer from the explicit expression of $W$ given by
\eqref{potential-W} that
$$
|W_{\eps}(v)| \leq |W(v)| + C v^2  ,  \quad |v|\leq \eps.
$$
On the other hand, a similar estimate holds in the
region $|v|\geq \eps$, thanks to the relation (\ref{relation-W}), so there exists a
positive $\eps$-independent constant $C$ such that
\begin{equation*}
|W_{\eps}(v)| \leq |W(v)| + C v^2  ,  \quad \text{for every }  v \in \mathbb{R}.
\end{equation*}
By Lebesgue's dominated convergence theorem, we have
$$
\int_{\mathbb{R}} W_{\eps}(v) dx \to \int_{\mathbb{R}} W(v) dx \quad
\text{as } \eps \to 0, \quad
\text{for every } v \in X.
$$
Therefore, $E_{\eps}(v_0) \to E(v_0)$ as $\eps \to 0$.

By using the above estimates and the bounds (\ref{bound-solution}), (\ref{bound-solution-2}), and (\ref{bound-solution-3}),
we obtain the following $\eps$- and $t$-uniform bound
for the sequence of approximating solutions $v^{\eps} \in C(\mathbb{R},H^1(\mathbb{R}))$
starting with the initial data $v_0 \in X$:
\begin{equation}
\label{bound-solution-eps}
\| v^{\eps}(t) \|_{H^1} + \| W_{\eps}(v^{\eps}(t)) \|_{L^1} \leq
C(P(v_0),E(v_0)) < \infty, \quad \forall t \in \mathbb{R}, \quad
\forall\eps \in (0,1],
\end{equation}
where the positive constant $C$ depends on the initial values $P(v_0)$ and $E(v_0)$ only. Therefore,
the sequence of approximating solutions $v^{\eps}$ is bounded in space
$L^{\infty}(\mathbb{R},X)$.
It also follows from the generalized KdV equation (\ref{generalizedKdV}) that the
sequence $v_t^{\eps}$ is bounded in space
$L^{\infty}(\mathbb{R},H^{-2}(\mathbb{R}))$. From Arzela--Ascoli Theorem, there exist
$v \in L^\infty(\mathbb{R},H^1(\mathbb R))$ and a subsequence of
$v^\eps \in L^\infty(\mathbb{R},H^1(\mathbb{R}))$, still denoted by $v^\eps$, such that
\begin{equation}
\label{limit-1}
v^\eps \to v \quad \text{strongly in } L^{\infty}_{\rm loc}(\mathbb{R},H^s_{\rm loc}(\mathbb R))
\quad  \mbox{\rm as} \;\; \eps \to 0, \quad {\rm for \; all} \;\; s < 1.
\end{equation}
Because the limit (\ref{limit-1}) includes the range
$\left(\frac{1}{2},1 \right)$ for $s$, up to subtracting another subsequence, we
may assume that
\begin{equation}
\label{limit-2}
v^{\eps}(x,t) \to v(x,t) \quad  \text{as }
\eps \to 0, \quad \text{for almost every } x \in \mathbb{R} \text{ and
  for every } t \in \mathbb{R}.
\end{equation}
In addition, because the upper bound (\ref{bound-solution-eps}) also
controls $\| W_{\eps}(v^{\eps}(t)) \|_{L^1}$, Fatou's lemma shows that
the limiting function $v$ belongs to $L^{\infty}(\mathbb{R},X)$.

To show that $P(v(t)) \leq P(v_0)$ and $E(v(t)) \leq E(v_0)$ for every
$t \in \mathbb{R}$, we use the weak lower semicontinuity of the $H^1$ norm and Fatou's lemma
for the potential term $[W_{\eps}]_-$ to obtain
$$
\| v(t) \|^2_{H^1} \leq \lim_{\eps \to 0} \| v^{\eps}(t) \|^2_{H^1}, \quad
\| [W_{\eps}(v(t))]_- \|_{L^1} \leq \lim_{\eps \to 0} \| [W_{\eps}(v^{\eps}(t))]_-  \|_{L^1},
\quad \mbox{\rm for \; every} \;\; t \in \mathbb{R}.
$$
On the other hand, using the relation (\ref{relation-W}), the limit (\ref{limit-2}),
and the fact that $W(v) \geq 0$ for $|v| \geq \sqrt{e}$, we obtain
$$
\| [W(v(t))]_+  \|_{L^1} =
\lim_{\eps \to 0} \| [W_{\eps}(v^{\eps}(t))]_+  \|_{L^1},
\quad \mbox{\rm for \; every} \;\; t \in \mathbb{R}.
$$
Then, it follows from the conservation
$$
P(v^{\eps}(t)) = P(v_0) \quad \mbox{\rm and} \quad
E_{\eps}(v^{\eps}(t)) = E_{\eps}(v_0)\quad \mbox{\rm for every} \;\; t \in \mathbb{R}
$$
and the convergence $E_{\eps}(v_0) \to E(v_0)$ as $\eps \to 0$ that
$$
P(v(t)) \leq P(v_0) \quad \mbox{\rm and} \quad
E(v(t)) \leq E(v_0)
$$

It remains to show that the limiting function $v \in L^{\infty}(\mathbb{R},X)$ is a weak solution of the log--KdV equation (\ref{logKdV}).
To do so, we first write a weak formulation of the generalized KdV equation (\ref{generalizedKdV}).
If $v^{\eps} \in C(\mathbb{R},H^1(\mathbb{R}))$ is a solution in Theorem \ref{theorem-KPV93} for $T = \infty$,
then for any test functions $\psi \in C^{\infty}_0(\mathbb{R}_x)$ and $\phi \in C^{\infty}_0(\mathbb{R}_t)$,
we have
\begin{equation}
\label{weak-generalized-KdV}
\int_{\mathbb{R}} \left[ \langle v^{\eps}, \psi \rangle_{L^2} \phi'(t) + \langle v^{\eps}, \psi''' \rangle_{L^2}
\phi(t) \right] dt + \int_{\mathbb{R}} \int_{\mathbb{R}} f_{\eps}(v^{\eps}) \psi'(x) \phi(t) dx dt = 0.
\end{equation}
Since $f_{\eps}(v) \to f(v)$ pointwise in $v$ as $\eps \to 0$, we apply the limit (\ref{limit-1}) to
the integral formulation (\ref{weak-generalized-KdV}) and obtain the following
integral equation for the limiting function $v$,
\begin{equation}
\label{weak-log-KdV}
\int_{\mathbb{R}} \left[ \langle v, \psi \rangle_{L^2} \phi'(t) + \langle v, \psi''' \rangle_{L^2}
\phi(t) \right] dt + \int_{\mathbb{R}} \int_{\mathbb{R}} f(v) \psi'(x) \phi(t) dx dt = 0,
\end{equation}
where $\psi$ and $\phi$ are any test functions. Therefore, $v \in L^{\infty}(\mathbb{R},X)$ is
a weak solution of the log--KdV equation (\ref{logKdV}), in particular, $v_t \in L^{\infty}(\mathbb{R},H^{-2}(\mathbb{R}))$.
The existence part of Theorem \ref{theorem-CP} is now proven.

\begin{rem}
\label{remark-L2}
Since $v_t$ is in $H^{-2}(\mathbb{R})$
(as opposed to the NLS case, where it is defined in
$H^{-1}(\mathbb{R})$, see \cite{Caz1,GV}),
we are not able to establish even conservation of the $L^2$ norm for the weak
solutions of the log--KdV equation (\ref{logKdV}), since we cannot prove rigorously
that
\begin{equation*}
  \frac{d}{dt}\|v(t)\|_{L^2}^2 = 2\langle v_t,v\rangle_{L^2} =0.
\end{equation*}
Nevertheless, we have proved that $\|v(t)\|_{L^2}^2$ is a non-increasing function of time $t$.
\end{rem}

\subsection{Uniqueness of solutions}

We will show uniqueness of the solution $v \in L^{\infty}(\mathbb{R},X)$
under the additional condition
$$
(\log|v|)_x \in L^{\infty}((-t_0,t_0) \times \mathbb{R}).
$$
Provided that uniqueness is proven, continuous dependence on the initial data $v_0 \in X$, the
$L^2$ norm and
energy conservation, and continuity of the solution $v \in C((-t_0,t_0),X)$ in time $t$ follow from
the arguments identical to the case of the log--NLS equation \cite[Section 9.3]{Caz}.

Assume existence of two local solutions $v$ and $u$ of the log--KdV equation (\ref{logKdV})
starting with the same initial data $v_0$.  Set $w := v-u$ such that $w|_{t = 0} = 0$.
Then $w$ satisfies a weak formulation similar to the integral equation (\ref{weak-log-KdV})
for the partial differential equation
$$
w_t + w_{xxx} + (v \log|v| - u \log|u|)_x = 0.
$$
Multiplying this equation by $w$, integrating over $x$, and
formally neglecting the values of $w$ and its derivatives as $|x| \to \infty$,
we obtain
$$
\frac{d}{dt} P(w) = -\int_{\mathbb{R}} (v_x \log|v| - u_x \log|u|) w dx,
$$
where $P(w) = \frac{1}{2} \| w\|^2_{L^2}$.
We use the following bound for the log-nonlinearity \cite{Caz}:
\begin{equation}
\label{log-bound}
|\log|v| - \log|u|| \leq \frac{|v-u|}{\min(|v|,|u|)}.
\end{equation}
We then write
\begin{align*}
  \frac{d}{dt} P(w) &= -\int_{|v|< |u|} (v_x \log|v| - u_x \log|u|)
  w dx - \int_{|v|\geq |u|} (v_x \log|v| - u_x \log|u|)
  w dx\\
&=  -\int_{|v|< |u|} v_x (\log|v| - \log|u|) w dx -
\int_{|v|< |u|} \log|u| w w_x dx\\
&\quad - \int_{|v|\geq |u|} u_x (\log|v| - \log|u|) w dx -
\int_{|v|\geq |u|} \log|v| w w_x dx.
\end{align*}
Applying (\ref{log-bound}), we obtain
\begin{align*}
 \left| \int_{|v|< |u|} v_x (\log|v| - \log|u|) w dx +\int_{|v|\geq
     |u|} u_x (\log|v| - \log|u|) w dx\right| \leq \left(\left\|
     \frac{v_x}{v} \right\|_{L^{\infty}} + \left\| \frac{u_x}{u}
   \right\|_{L^{\infty}} \right)2P(w).
\end{align*}
Integrating by parts, we also have
\begin{equation*}
-\int_{|v|< |u|} \log|u| w w_x dx - \int_{|v|\geq |u|} \log|v| w w_x dx =
\frac{1}{2} \int_{|v|< |u|} (\log|u|)_x w^2 dx +
\frac{1}{2} \int_{|v|\geq |u|} (\log|v|)_x w^2 dx,
\end{equation*}
hence
\begin{equation*}
\left| \frac{d}{d t}  P(w) \right| \leq
3\left(  \left\| \frac{v_x}{v} \right\|_{L^{\infty}} + \left\| \frac{u_x}{u} \right\|_{L^{\infty}} \right) P(w).
\end{equation*}
Gronwall's inequality implies that if $w |_{t = 0} = 0$, then $P(w) = 0$
for all $t \in (-t_0,t_0)$, for which the solutions $v,u$ of the log--KdV equation (\ref{logKdV})
satisfy $(\log|v|)_x, (\log |u|)_x \in L^{\infty}((-t_0,t_0)\times\mathbb{R})$.
The uniqueness part of Theorem \ref{theorem-CP} is now proven.

\section{Spectral stability of Gaussian solitary waves}

To prove Theorem \ref{theorem-spectral}, we study the spectrum of the linearized operator
$A := \partial_x L$ given by
\begin{equation}
A = -\partial_x^3 + \frac{1}{4} (x^2 - 6) \partial_x + \frac{1}{2} x.
\end{equation}
The domain of this operator with a range in $L^2(\mathbb{R})$ is given by
\begin{equation}
D(A) = \left\{ u \in H^3(\mathbb{R}) : \quad x^2 \partial_x  u \in L^2(\mathbb{R}), \quad x u \in L^2(\mathbb{R}) \right\}.
\end{equation}
We note that the eigenfunctions of the linearized operator $\partial_x L$ associated with the
KdV-type spectral problem $A f = \lambda f$ is defined in the function space $X_A := D(A) \cap \dot{H}^{-1}(\mathbb{R})$,
hence the anti-derivative of the eigenfunction $f$ is required to be squared integrable \cite{kapstef,dmitrystab}.

We shall employ the Fourier transform $\mathcal{F} : L^2(\mathbb{R}) \to L^2(\mathbb{R})$ defined by
\begin{equation}
\label{Fourier}
\hat{u}(k) := \mathcal{F}(u)(k) = \frac{1}{\sqrt{2\pi}} \int_{\mathbb{R}} u(x) e^{-ikx} dx, \quad k \in \mathbb{R}.
\end{equation}
The Fourier transform $\mathcal{F}$ is helpful in the study of the spectrum of $A$ because of dualism between
derivatives and multiplications. After an application of the Fourier transform $\mathcal{F}$, the third-order
differential operator $A$ in $x$-space is mapped to the second-order differential operator $\hat{A}$ in $k$-space,
where
\begin{equation}
\hat{A} = \frac{i}{4} k \left( -\partial_k^2 + 4 k^2 - 6 \right).
\end{equation}
The function space $X_A = D(A) \cap \dot{H}^{-1}(\mathbb{R})$ is now mapped to the function space
$\hat{X}_A$ given by
\begin{equation}
\label{function-space-hat-X}
\hat{X}_A = \left\{ \hat{u} \in H^1(\mathbb{R}) : \quad k \partial_k^2 \hat{u} \in L^2(\mathbb{R}), \quad
k^3 \hat{u} \in L^2(\mathbb{R}), \quad k^{-1} \hat{u} \in L^2(\mathbb{R}) \right\}.
\end{equation}
Denoting $\hat{B} := k (-\partial_k^2 + 4 k^2 - 6)$, we have $\sigma(\hat{A}) = \frac{i}{4} \sigma(\hat{B})$.

\subsection{Double zero eigenvalue of $\hat{A}$}

Let $\hat{v}_G (k):= e^{-k^2}$ be the Fourier transform of the Gaussian solitary wave $v_G$
(up to the constant multiplicative factor).
We check by direct computation that
\begin{equation}
\label{double-zero}
\hat{B} \partial_k \hat{v}_G = 0 \quad {\rm and} \quad \hat{B} \hat{v}_G = - 4k \hat{v}_G = 2 \partial_k \hat{v}_G.
\end{equation}
Since $\hat{B} = 4 k \hat{L}$, where $\hat{L}$ is the Fourier image of the Schr\"{o}dinger operator
$L$ with a harmonic potential (\ref{Schrodinger-operator}),  we conclude that
${\rm null}(\hat{B}) = {\rm span}\{\partial_k \hat{v}_G \}$, whereas the generalized null space of $\hat{B}$
includes the two-dimensional subspace $\hat{X}_0 = {\rm span}\{\partial_k \hat{v}_G,\hat{v}_G \}$.

Moreover, since $\hat{B}^* \hat{v}_G = 0$ and $\| \hat{v}_G \|_{L^2} \neq 0$,
no solution $\hat{u} \in \hat{X}_A$ of the inhomogeneous equation $\hat{B} \hat{u} = \hat{v}_G$
exists. Therefore, the Jordan chain for the zero eigenvalue of $\hat{B}$ is two-dimensional and
the generalized null space of $\hat{B}$ is exactly $\hat{X}_0$.
With the inverse Fourier transform, these arguments conclude
consideration of the zero eigenvalue in Theorem \ref{theorem-spectral}.

\subsection{Nonzero eigenvalues of $\hat{A}$}

Let us now consider nonzero values for the spectral parameter $\lambda$ in the spectral
problem $\hat{A} \hat{f} = \lambda \hat{f}$. Recall that
$\sigma(\hat{A}) = \frac{i}{4} \sigma(\hat{B})$. Therefore, let $\lambda = \frac{i}{4} E$
and consider the spectral problem $\hat{B} \hat{u} = E \hat{u}$ rewritten as
the following differential equation
\begin{equation}
\label{spectral-problem}
\frac{d^2 \hat{u}}{d k^2} + \left( \frac{E}{k} + 6 - 4 k^2 \right) \hat{u}(k) = 0, \quad k \in \mathbb{R}.
\end{equation}
We employ the theory of differential equations to study solutions of the spectral problem
(\ref{spectral-problem}).

If $E \neq 0$, the point $k = 0$ is a regular singular point
of the differential equation (\ref{spectral-problem}) with indices $0$ and $1$.
By Frobenius' method \cite[Chapter 4]{Teschl},
there exist two linearly independent solutions of this differential equation.
The first solution is given by the power series expansion
\begin{equation}
\label{expansion-1}
\hat{u}_1(k) = k -\frac{1}{2} E k^2 + \left( \frac{1}{12} E^2 - 1 \right) k^3 + \mathcal{O}(k^4) \quad \mbox{\rm as} \quad k \to 0
\end{equation}
and the other solution is given by the logarithm-modified power series expansion
\begin{equation}
\label{expansion-2}
\hat{u}_2(k) = 1 - 3 \left(1 + \frac{1}{4} E^2 \right) k^2 + \mathcal{O}(k^3)
- E k \log(k) \left( 1 - \frac{1}{2} E k + \mathcal{O}(k^2) \right) \quad \mbox{\rm as} \quad k \to 0.
\end{equation}
Because the other singular point of the differential equation (\ref{spectral-problem}) is infinity,
the expansions (\ref{expansion-1}) and (\ref{expansion-2}) converge for every $k \in \mathbb{R}$.

Recall again the function space $\hat{X}_A$ defined by (\ref{function-space-hat-X}).
Because $k^{-1} \hat{u} \in L^2(\mathbb{R})$ is required in $\hat{X}_k$,
we have $\hat{u} \in \hat{X}_A$ if and only if $\hat{u}$ is constant proportional to the solution $\hat{u}_1$.
Hence, both for $k > 0$ and $k < 0$, the solution of the differential equation (\ref{spectral-problem}) is uniquely defined by the
solution $\hat{u}_1$ (up to the constant multiplicative factor) for every $E \in \mathbb{C} \backslash \{0\}$.
The admissible values of $E$ are determined from the behavior of the solution $\hat{u}(k)$ as $|k| \to \infty$.

By the WKB method without turning points \cite[Chapter 7.2]{Miller}, there exist two linearly independent
solutions of the differential equation (\ref{spectral-problem}): one solution diverges and the other one
decays to zero as $|k| \to \infty$. The decaying solution is defined by the asymptotic behavior
\begin{equation}
\label{decay-behavior}
\hat{u}(k) \sim k e^{-k^2} \left( 1 - \frac{E}{4k} + \mathcal{O}\left(\frac{1}{k^2}\right) \right) \quad
\mbox{\rm as} \quad |k| \to \infty.
\end{equation}
It becomes now clear that two decay conditions at $k \to +\infty$ and $k \to -\infty$ over-determine
the spectral problem (\ref{spectral-problem}) because the behavior of $\hat{u}$
both for $k > 0$ and $k < 0$ is uniquely determined by only one spectral parameter $E$
(up to the constant multiplicative factors) and no parity symmetry exists for $k > 0$ and $k < 0$.

The way around this obstacle is to consider a weak formulation for the solutions of the
differential equation (\ref{spectral-problem}) in the function space $\hat{X}_A$. Since the
differential equation with $E \neq 0$ has singularities at $\{-\infty,0,\infty\}$, we split $\mathbb{R}$ into two sets
$I_- := (-\infty,0)$ and $I_+ := (0,\infty)$ and look for piecewisely defined eigenfunctions
$\hat{u}_{\pm}$ supported on $I_{\pm}$ only.

Let $\hat{B}_{\pm}$ denote the operator $\hat{B}$ restricted on $I_{\pm}$ subject to the Dirichlet boundary condition at
$k = 0$. By Theorem 4 in \cite[p.1438]{DS}, the continuous spectrum of $B$ is the union of
the continuous spectra of operators $B_+$ and $B_-$.

We shall first characterize the spectrum of the operator $\hat{B}_+$ on $I_+$ and prove that
this spectrum is purely discrete. By using a substitution
$\hat{u}_+(k) = k^{1/2} \hat{v}_+(k)$ for $k > 0$, we reduce the spectral problem $\hat{B}_+ \hat{u}_+ = E \hat{u}_+$
to the symmetric form
\begin{equation}
\label{spectral-problem-plus}
k^{1/2} \left( - \frac{d^2}{d k^2} + 4 k^2 - 6 \right) k^{1/2} \hat{v}_+(k) = E \hat{v}_+(k), \quad k \in I_+.
\end{equation}
Since $k = 0$ is singular, we require $|\hat{v}_+(0)| < \infty$, although the solution $\hat{u}_1$
in (\ref{expansion-1}) implies more precisely that $\hat{v}_+(k) = \mathcal{O}(k^{1/2})$ as $k \to 0$.
The spectral problem (\ref{spectral-problem-plus}) in the function space with $|\hat{v}_+(0)| < \infty$
is self-adjoint, hence the admissible values of $E$ are real.

We shall now claim that $\sigma(\hat{B}_+) = \{ E_n \}_{n \in \mathbb{N}_+}$, where
$0 = E_0 < E_1 < E_2 < ... $ and $E_n \to +\infty$ as $n \to \infty$.
First, because the potential $4 k^2 - 6$ is confining, the resolvent operator is compact,
and the spectrum of $\hat{B}_+$ is purely discrete.
Next, $\hat{v}_+ = k^{1/2} e^{-k^2}$ is the exact solution of
the spectral problem (\ref{spectral-problem-plus}) for $E = 0$
and this eigenfunction is strictly positive for all $k \in I_+$. By Sturm's Theorem,
$E_0 = 0$ is at the bottom of the spectrum of $\hat{B}_+$, which is bounded from below,
and hence $E_n > 0$ for every $n \in \mathbb{N}$. Finally, the eigenfunctions are
uniquely determined by the expansion (\ref{expansion-1}) as $k \to 0$, hence
each eigenvalue is simple.

As a consequence of the reduction to the self-adjoint problem (\ref{spectral-problem-plus}),
we have orthogonality and normalization of the real-valued eigenfunctions $\hat{u}_{+n}(k) = k^{1/2} \hat{v}_{+n}(k)$
corresponding to the real eigenvalues $E_n$ for all $n \in \mathbb{N}_+$:
\begin{equation}
\label{orthogonality-normalization-plus}
\delta_{n,m} = \int_0^{\infty} \hat{v}_{+n}(k) \hat{v}_{+m}(k) dk = \int_0^{\infty} \frac{\hat{u}_{+ n}(k) \hat{u}_{+m}(k)}{k} dk,
\end{equation}
where $\delta_{n,m}$ is the Kronecker's symbol.

The same analysis applies to the spectrum of the operator $\hat{B}_-$ on $I_-$
with the only difference that the substitution $\hat{u}_-(k) = (-k)^{1/2} \hat{v}_-(k)$ for $k < 0$
reduces the spectral problem $\hat{B}_- \hat{u}_- = E \hat{u}_-$
to the symmetric form
\begin{equation}
\label{spectral-problem-minus}
(-k)^{1/2} \left( - \frac{d^2}{d k^2} + 4 k^2 - 6 \right) (-k)^{1/2} \hat{v}_-(k) = -E \hat{v}_-(k), \quad k \in I_-,
\end{equation}
subject to the boundary behavior $|\hat{v}_-(0)| < \infty$.
Therefore, $\sigma(\hat{B}_-) = \{ -E_n \}_{n \in \mathbb{N}_+}$, where
$E_n$ are the same eigenvalues as in the spectral problem (\ref{spectral-problem-plus}).
Similarly to (\ref{orthogonality-normalization-plus}), it follows from the self-adjoint problem (\ref{spectral-problem-minus})
that the real-valued eigenfunctions $\hat{u}_{-n}(k) = (-k)^{1/2} \hat{v}_{-n}(k)$
corresponding to the real eigenvalues $-E_n$ for all $n \in \mathbb{N}_+$ are orthogonal and can be normalized by
\begin{equation}
\label{orthogonality-normalization-minus}
\delta_{n,m} = \int_{-\infty}^0 \hat{v}_{-n}(k) \hat{v}_{-m}(k) dk = -\int_{-\infty}^0 \frac{\hat{u}_{- n}(k) \hat{u}_{-m}(k)}{k} dk.
\end{equation}

The statement about nonzero eigenvalues of the operator $A$ in Theorem \ref{theorem-spectral} is now proven with
the correspondence $\omega_n = \frac{1}{4} E_n$ for all $n \in \mathbb{N}$.

\begin{rem}
The double zero eigenvalue associated with the two-dimensional subspace
$\hat{X}_0 = {\rm span}\{\partial_k \hat{v}_G,\hat{v}_G \}$ for the operator $\hat{A}$
is mapped into a semi-simple zero eigenvalue associated with the simple
zero eigenvalues of the operators $\hat{B}_+$ and $\hat{B}_-$. There is no
ambiguity here because if $E = 0$, then $k = 0$ is an ordinary point
of the differential equation (\ref{spectral-problem}) and the splitting of $\mathbb{R}$
into $I_-$ and $I_+$ does not make sense.
\end{rem}

\subsection{Eigenfunctions of the spectral problem}

We now confirm that the eigenfunction $\hat{u}$ defined piecewise as
$$
\mbox{\rm either} \quad \hat{u}(k) = \left\{ \begin{array}{l}\hat{u}_+(k), \quad k > 0, \\
0, \quad \quad \quad k < 0, \end{array} \right. \quad \mbox{\rm or} \quad
\hat{u}(k) = \left\{ \begin{array}{l} 0, \quad \quad \quad k > 0, \\
\hat{u}_-(k), \quad k < 0, \end{array} \right.
$$
belongs to the space $\hat{X}_A$. We note that $\hat{u}$ is continuous on $\mathbb{R}$
and piecewise $C^1$ with the jump discontinuity of the first derivative at $k = 0$.
Moreover, $\hat{u}(k) = \mathcal{O}(k)$ as $k \to 0$. Therefore,
$\hat{u} \in H^1(\mathbb{R})$ and $k^{-1} \hat{u} \in L^2(\mathbb{R})$.
Furthermore, $\hat{u}(k)$ decays to zero as $|k| \to \infty$ according to
the asymptotic behavior (\ref{decay-behavior}). Therefore, $k^3 \hat{u} \in L^2(\mathbb{R})$.
Finally, $\partial_k^2 \hat{u}$ is proportional to the Dirac delta distribution $\delta$ at $k = 0$
and is smooth for $k \neq 0$, hence
$\partial_k^2 \hat{u} \notin L^2(\mathbb{R})$ but $k \partial_k^2 \hat{u} \in L^2(\mathbb{R})$.
Thus, $\hat{u} \in \hat{X}_A$.

The Fourier transform $\mathcal{F}$ is dual between the smoothness and decay of the eigenfunctions.
Since the eigenfunctions $\hat{u}(k)$ decays fast as $|k| \to \infty$, according to
the Gaussian decay (\ref{decay-behavior}), their inverse Fourier transform
$u(x)$ is smooth in $x$. Moreover, $u_{\pm}$ defined by
\begin{equation}
\label{inverse-Fourier}
u_+(x) = \frac{1}{\sqrt{2\pi}} \int_0^{\infty} \hat{u}_+(k) e^{i k x} dk, \quad
u_-(x) = \frac{1}{\sqrt{2\pi}} \int_{-\infty}^0 \hat{u}_-(k) e^{i k x} dk
\end{equation}
are analytically extended for ${\rm Im}(x) \gtrless 0$ respectively. On
the other hand, $\hat{u}(k)$ have jump discontinuity in the first derivative across $k = 0$,
therefore, $u(x)$ decay only algebraically as $|x| \to \infty$. The only exception
is the double zero eigenvalue, since one can glue $\hat{u}_+(k)$ and $\hat{u}_-(k)$
for the eigenvalue $E_0 = 0$ to a smooth eigenfunction $\hat{u}(k) = k e^{-k^2}$,
which corresponds to the eigenfunction $u(x) = x e^{-x^2/4}$ with the fast decay as $|x| \to \infty$.

The statement about eigenfunctions of the operator $A$ in Theorem \ref{theorem-spectral} is now proven.

\subsection{Numerical illustrations}

\begin{figure}[h]
\begin{center}
\includegraphics[scale=0.7]{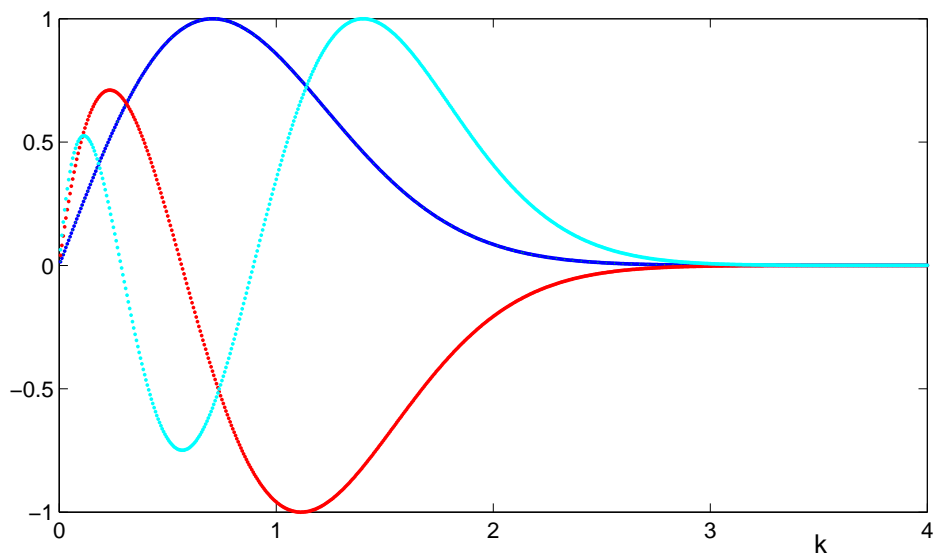}
\end{center}
\caption{\label{fig-spectrum}
Eigenfunctions $\hat{u}$ of the spectral problem (\ref{spectral-problem}) versus $k$
for the first three eigenvalues $0 = E_0 < E_1 < E_2$. Eigenfunctions for $E_1$ and $E_2$ have
one and two zeros for $k > 0$, respectively.}
\end{figure}

Figure \ref{fig-spectrum}  shows the first three eigenfunctions $\hat{u}_+(k) = k^{1/2} \hat{v}_+(k)$ of the spectral problem
(\ref{spectral-problem-plus}) for the first three eigenvalues $E_0 = 0$, $E_1 \approx 5.4109$,
and $E_2 \approx 12.3080$. These eigenfunctions and eigenvalues were computed
numerically by means of the central difference approximation of the derivatives
in the spectral problem (\ref{spectral-problem-plus}) and the MATLAB eigenvalue solver.
The fast decay of the eigenfunctions $\hat{u}_+(k)$ as $k \to \infty$ and the Sturm's nodal properties of
the eigenfunctions for $k > 0$ are obvious from the figure.

\begin{figure}[h]
\begin{center}
\includegraphics[scale=0.37]{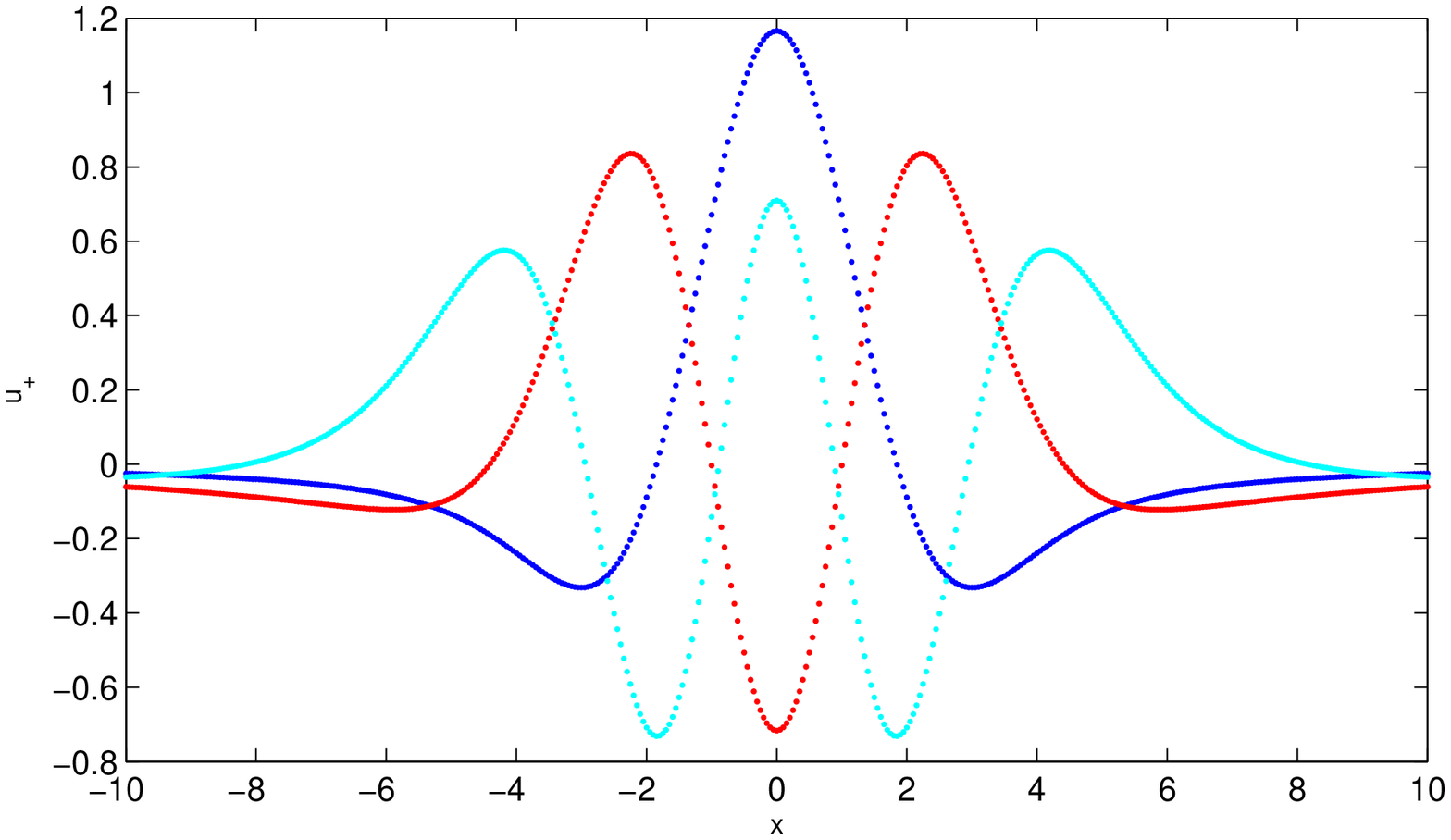}
\includegraphics[scale=0.37]{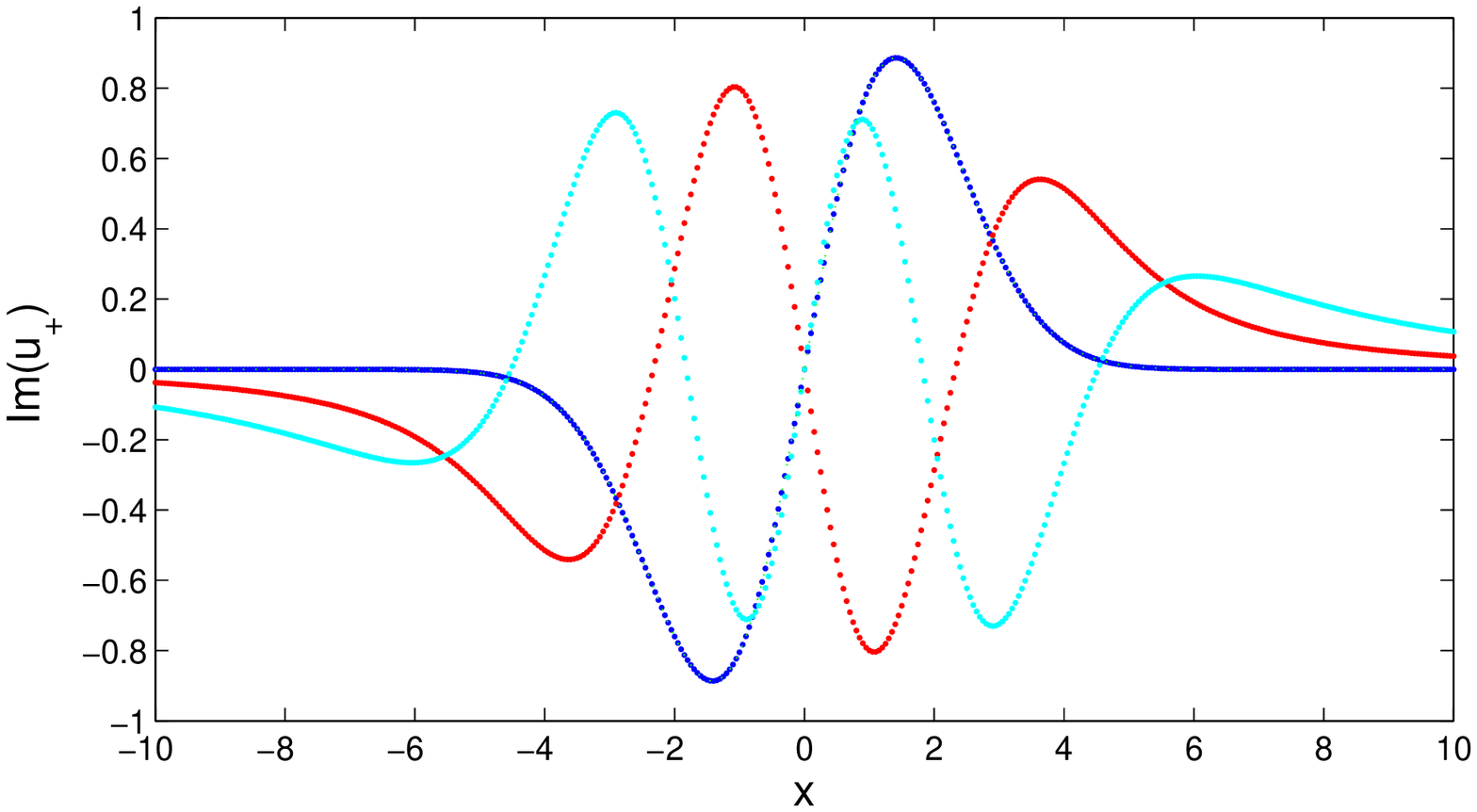}
\end{center}
\caption{\label{fig-functions}
Real (left) and imaginary (right) parts of the three eigenfunctions $u$ versus $x$
for the first three eigenvalues $E$.}
\end{figure}

Figure \ref{fig-functions} shows the real (left) and imaginary (right)
parts of the eigenfunctions $u_+$ versus $x$ after the inverse Fourier transform (\ref{inverse-Fourier}).
The imaginary part of the first eigenfunction decays fast as $|x| \to \infty$, according to the exact
eigenfunction $v_G'(x) \sim x e^{-x^2/4}$, which is also
shown by a dotted line (invisible from the numerical dots).

To inspect the slow (algebraic) decay of the eigenfunctions, we multiply
the eigenfunctions by the factor $(1+x^2)^{p/2}$ with $p = 2$ and $p = 3$.
It follows from the numerical data that the real
parts of the eigenfunctions decay like $\mathcal{O}(|x|^{-2})$ as $|x| \to \infty$,
whereas the imaginary parts of the eigenfunctions (except for the first eigenfunction)
decay like $\mathcal{O}(|x|^{-3})$ as $|x| \to \infty$. This corresponds
to the finite jump discontinuity of the first derivative of $\hat{u}(k)$ at $k = 0$.

\section{Cauchy problem for the linearized log--KdV equation}

We consider the Cauchy problem for the linearized log-KdV equation (\ref{linlogKdV})
at the Gaussian solitary wave $v_G$. For convenience, we rewrite the Cauchy problem again:
\begin{equation}
\label{linlogKdV-section}
\left\{ \begin{array}{l} u_t = \partial_x L u, \quad t > 0, \\
u |_{t = 0} = u_0, \end{array} \right.
\end{equation}
We shall first make use of Theorem \ref{theorem-spectral} to give a quick proof of Theorem \ref{theorem-GP}.
This proof is alternative to the arguments given in \cite{JP13}, which rely on the symplectic projections
and the energy method. Next, we shall approximate the Cauchy problem
numerically, starting with the Gaussian initial data. We conclude the paper
with a short discussion.

\subsection{Proof of Theorem \ref{theorem-GP}}
\label{sec:proofGP}

Let $\{ u_{\pm n}(x) \}_{n \in \mathbb{N}}$ be a sequence of eigenfunctions of $\partial_x L$
for the sequence of nonzero imaginary eigenvalues $\{ \pm i \omega_n \}_{n \in \mathbb{N}}$
constructed in Theorem \ref{theorem-spectral}. It follows from the orthogonality and
normalization conditions (\ref{orthogonality-normalization-plus}) and (\ref{orthogonality-normalization-minus})
that the eigenfunctions are orthogonal to each other
and normalized with respect to the symplectic inner product:
\begin{equation}
\label{symplectic-orthogonality}
\langle \partial_{x}^{-1} u_{\pm n}, u_{\pm m} \rangle_{L^2} :=
\int_{-\infty}^{\infty} \frac{\hat{\bar{u}}_{\pm n}(k) \hat{u}_{\pm m}(k)}{i k} dk = \mp i \delta_{n,m},
\end{equation}
where $\partial_x^{-1} u(x) := \int_{-\infty}^{x} u(x') dx'$
and the Fourier transform (\ref{Fourier}) is used. Note that
$u_{\pm n} \in L^1(\mathbb{R})$ if $u_{\pm n} \in L^2_1(\mathbb{R}) \subset D(A)$.
Also note that $\hat{u}_{\pm n}(0) = 0$ if $u_{\pm n} \in \dot{H}^{-1}(\mathbb{R}) \subset X_A$.

Let us consider the decomposition of the solution $u$ of the linearized log--KdV equation (\ref{linlogKdV-section})
as the series of eigenfunctions of the operator $A$:
\begin{equation}
\label{decomposition}
u(x,t) = b \left[ v_G(x) - t v_G'(x) \right] + a_0 v_G'(x) + \sum_{n \in \mathbb{N}} a_{+n} u_{+n}(x) e^{i \omega_n t} +
\sum_{n \in \mathbb{N}} a_{-n} u_{-n}(x) e^{-i \omega_n t},
\end{equation}
where the coefficients $b$ and $\{ a_n \}_{n \in \mathbb{Z}}$ are found uniquely from the
conditions of symplectic orthogonality and normalization (\ref{symplectic-orthogonality})
applied to the initial data $u |_{t = 0} = u_0$:
\begin{equation}
\label{projection-1}
b = \frac{\langle v_G, u_0 \rangle_{L^2}}{\langle v_G, v_G \rangle_{L^2}}, \quad
a_0 = \frac{1}{\langle v_G, v_G \rangle_{L^2}} \left[ \frac{b}{2}
\left(\int_{-\infty}^{\infty} v_G(x) dx \right)^2 - \langle \partial_x^{-1} v_G, u_0 \rangle_{L^2} \right),
\end{equation}
and
\begin{equation}
\label{projection-2}
a_{\pm n} = \pm i \langle \partial_x^{-1} u_{\pm n}, u_0 \rangle_{L^2}, \quad n \in \mathbb{N}.
\end{equation}
If $u_0 \in X_c$ defined in (\ref{constrained-space}), then $b = 0$.
The coefficients $\{ a_n \}_{n \in \mathbb{Z}}$ are well defined for $u_0 \in X_c$.
Moreover, by the spectral theorem and the reduction of the non-self-adjoint operator $A$
to the two self-adjoint operators in Section 3.2, the series (\ref{decomposition})
converges in the $L^2$ sense for every $t \in \mathbb{R}$.

Define now the conserved energy of the linearized log--KdV equation,
\begin{equation}
\label{projection-3}
E_c(u) := \frac{1}{2} \langle L u, u \rangle_{L^2} = {\rm const}, \quad t \in \mathbb{R}.
\end{equation}
Recall that $E_c(u_0) \geq 0$ if $u_0 \in X_c$.
It follows from (\ref{decomposition}) and (\ref{projection-3}) that
if $u_0 \in X_c$, then $u(\cdot,t) \in X_c$ and $E_c(u(\cdot,t)) \geq 0$
for every $t \in \mathbb{R}$. Let $y(x,t) := u(x,t) - a_0 v_G'(x)$
and recall that $b = 0$ for $u_0 \in X_c$. Under these conditions, we have
$$
E_c(y(\cdot,t)) = E_c(u(\cdot,t)) = E_c(u_0), \quad t \in \mathbb{R}
$$
and
$$
E_c(y(\cdot,t)) = \frac{1}{2} \sum_{n \in \mathbb{N}} \omega_n \left( |a_{+,n}|^2 + |a_{-n}|^2 \right),
$$
where the orthogonality and normalization conditions (\ref{symplectic-orthogonality}) have been used.

Because $y$ belongs to a subspace of $X_c$ spanned by $\{ u_{\pm n} \}_{n \in \mathbb{N}}$
associated with the nonzero eigenvalues $\{ \pm i \omega_n \}_{n \in \mathbb{N}}$, where each $\omega_n > 0$,
the energy functional $E_c(y)$ is equivalent to the squared $H^1(\mathbb{R}) \cap L^2_1(\mathbb{R})$ norm.
Therefore, there is a $t$-independent positive constant $C$ such that
\begin{equation}
\label{projection-4}
\| y(\cdot,t) \|^2_{H^1} + \| y(\cdot,t) \|_{L^2_1}^2 \leq C E_c(y(\cdot,t)) = C E_c(u_0) \leq
C \left( \frac{1}{2} \| u_0 \|_{H^1}^2 + \frac{1}{8} \| u_0 \|_{L^2_1}^2 \right).
\end{equation}
By the decomposition (\ref{decomposition}), the bound (\ref{projection-4}) and the triangle inequality,
we obtain the bound (\ref{global-bound}) for every $t \in \mathbb{R}$. Theorem \ref{theorem-GP} is now proven.

\subsection{Numerical illustrations}

We truncate the spatial domain on $[-L,L]$ for sufficiently large $L$ ($L = 40$ was used
in our numerical results) and discretize $u(x,t)$ at equally spaced grid points $\{ x_k \}_{k = 0}^{N+1}$
subject to the periodic boundary conditions $u(x_0,t) = u(x_{N+1},t)$. Using the central
difference approximation and the Heun's method, we rewrite the linearized log--KdV equation
(\ref{linlogKdV}) in the iterative form
\begin{equation}
\label{numerical-scheme}
\left( I - \frac{\Delta t}{2} D L_D \right) {\bf u}_{m+1} = \left( I + \frac{\Delta t}{2} D L_D \right) {\bf u}_m
\end{equation}
where $\Delta t$ is the time step, ${\bf u}_m$ is the vector of discretized solution at the time $t_m$,
$I$ is an identity matrix,
$D$ is the matrix for the first derivative, and $L_D$ is the matrix for the second-order differential operator $L$.

For a positive parameter $\alpha$, we consider the odd initial data
\begin{equation}
\label{odd-data}
u_0 (x)= x \left(x^2 - \frac{5+6\alpha}{(1+\alpha)(1+2\alpha)}\right)
e^{-\frac{(1+2\alpha) x^2}{4}},
\end{equation}
which satisfies the symplectic orthogonality conditions
$\langle v_G, u_0 \rangle_{L^2} = \langle \partial_x^{-1}v_G, u_0 \rangle_{L^2} = 0$.
These constraints ensure that $b = 0$ and $a_0 = 0$ in the decomposition (\ref{decomposition})
and the solution $u$ of the linearized log--KdV equation (\ref{linlogKdV-section})
is spanned by eigenfunctions of the linearized operator $A = \partial_x L$ for
nonzero eigenvalues.

\begin{figure}[h]
\begin{center}
\includegraphics[scale=0.35]{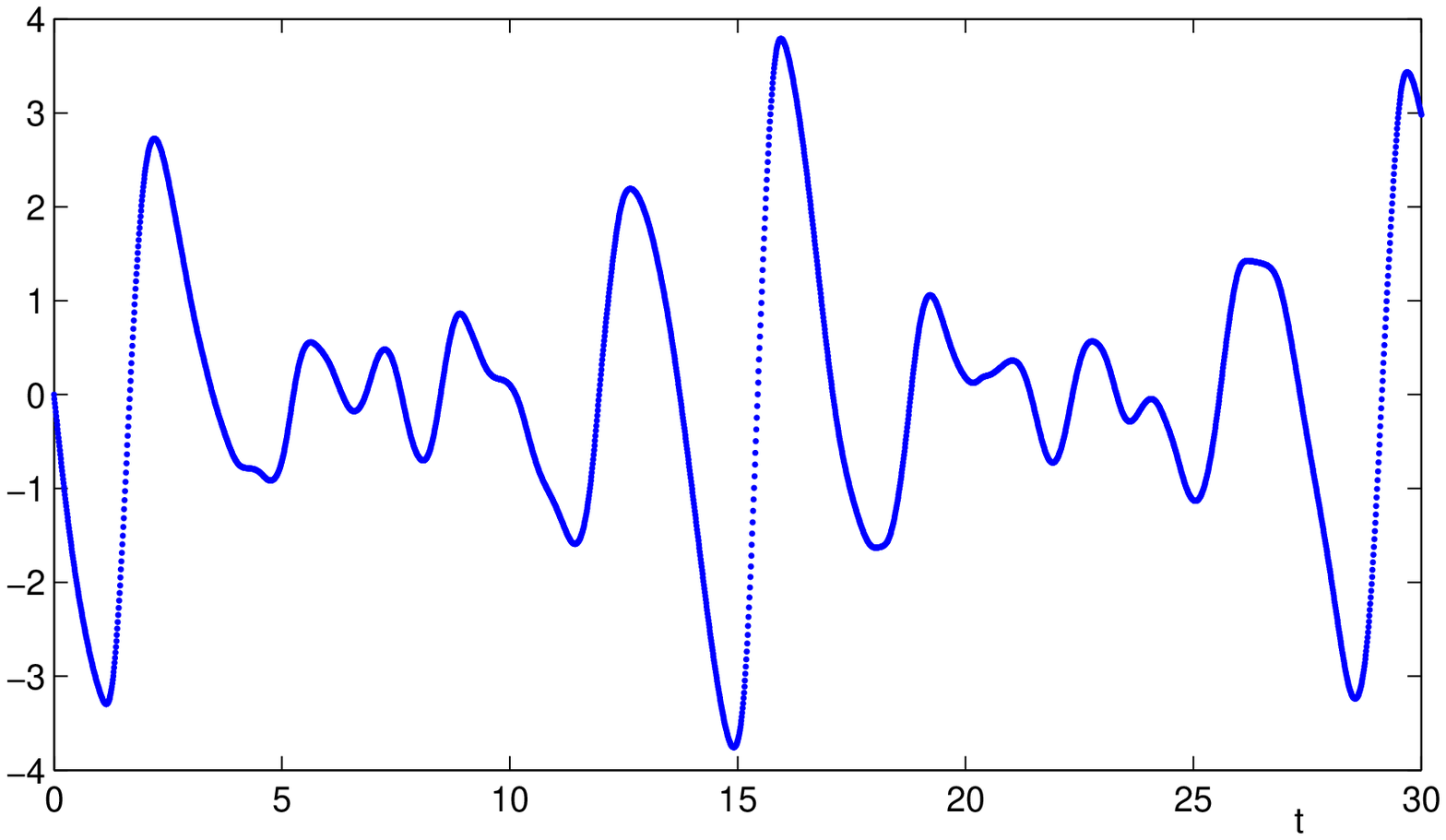}
\includegraphics[scale=0.35]{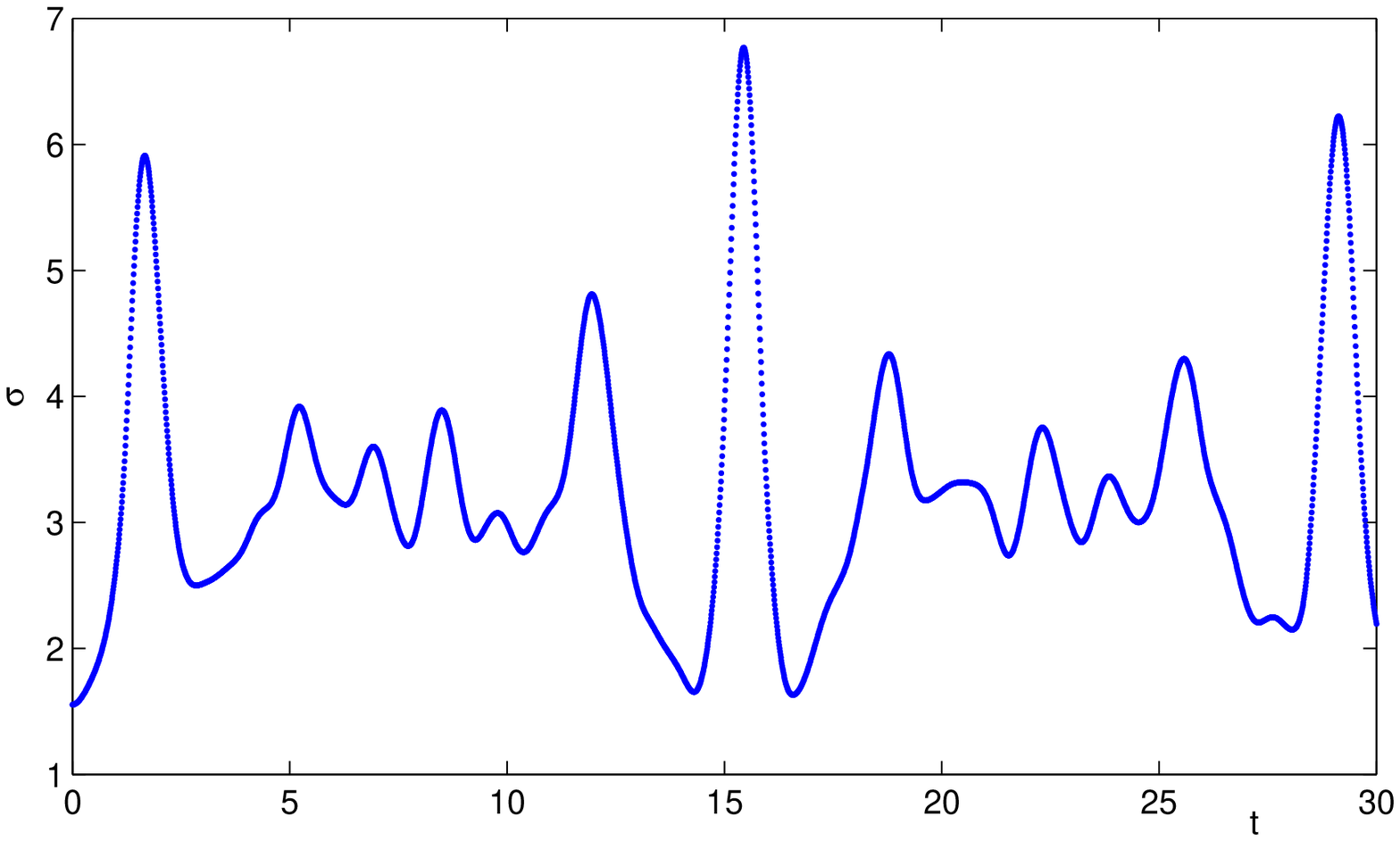}
\includegraphics[scale=0.4]{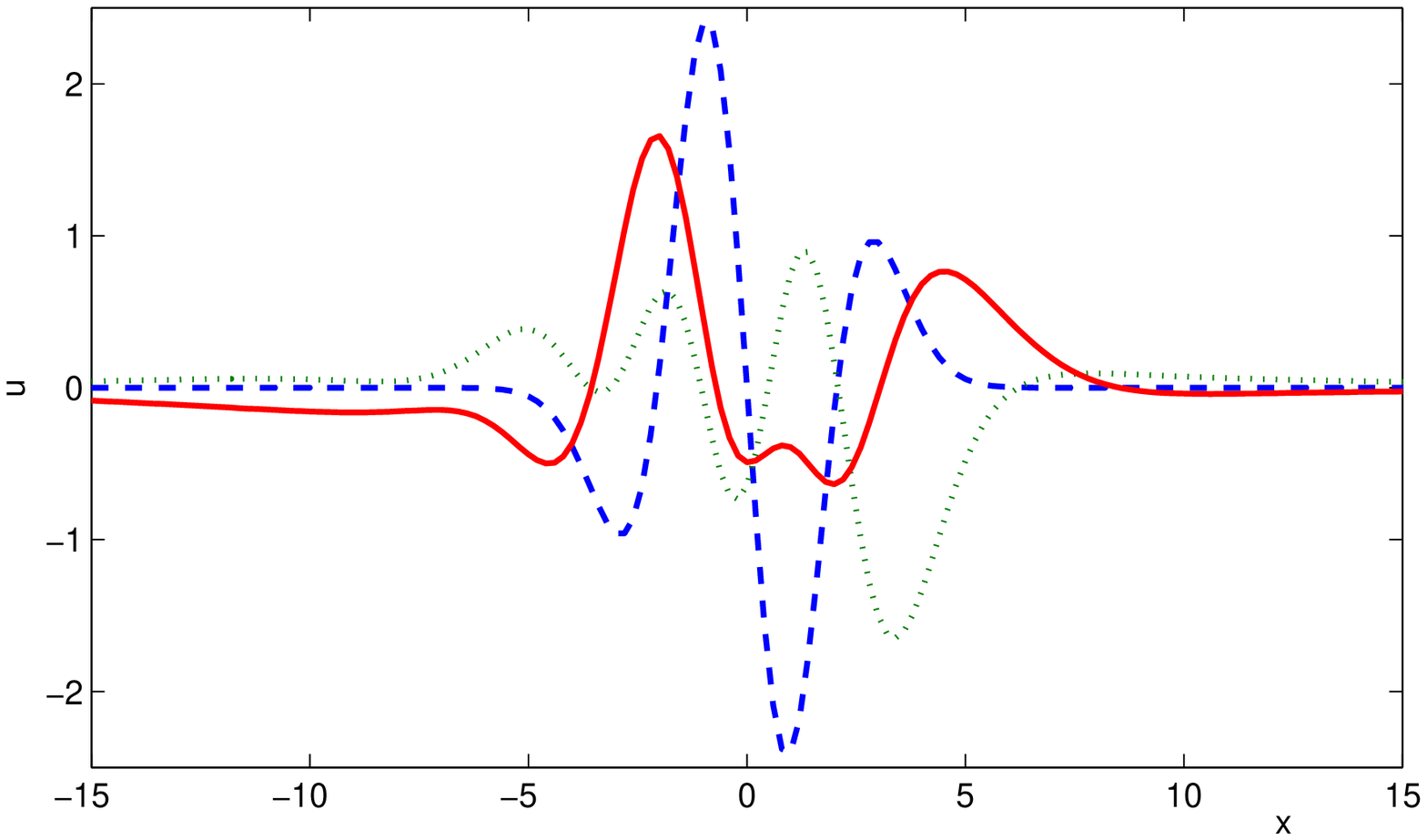}
\end{center}
\caption{\label{fig-Cauchy-2}
Numerical solution of the linearized log--KdV equation
(\ref{linlogKdV}) with an odd initial data (\ref{odd-data}):
center of mass $\bar{x}$ (top left) and standard deviation $\sigma$ (top right) versus time;
the profile $u(x,t)$ versus $x$ (bottom) for
$t = 0$ (blue dashed), $t = 2.5$ (green dotted), and $t = 5$ (red solid).}
\end{figure}

Figure \ref{fig-Cauchy-2} reports numerical computations of the
iterative scheme (\ref{numerical-scheme})
for the Gaussian initial data (\ref{odd-data}) with $\alpha = 0.1$. Besides the profiles
of the solution shown on the spatial interval $[-15,15]$ (bottom panel), we
also show the center of mass (top left panel)
\begin{equation*}
  \bar{x}(t) = \frac{\int_{\mathbb{R}} x u^2(x,t) dx}{\int_{\mathbb{R}} u^2(x,t) dx},
\end{equation*}
and the standard deviation (top right panel)
$$
\sigma(t) = \left( \frac{\int_{\mathbb{R}} (x - \bar{x}(t))^2 u^2(x,t) dx}{\int_{\mathbb{R}} u^2(x,t) dx} \right)^{1/2}.
$$
We can see from the top right panel of Figure \ref{fig-Cauchy-2} that the standard deviation $\sigma(t)$ oscillates periodically
with a large amplitude, which is still much smaller than the half-size of the computational domain $L = 40$.
Therefore, the odd Gaussian data (\ref{odd-data}) evolves into a solution, which does not spread
out in the time evolution of the linearized  log--KdV equation
(\ref{linlogKdV}). Nevertheless, a visible radiation appears on the left slope of the Gaussian pulse.

\begin{figure}[h]
\begin{center}
\includegraphics[scale=0.35]{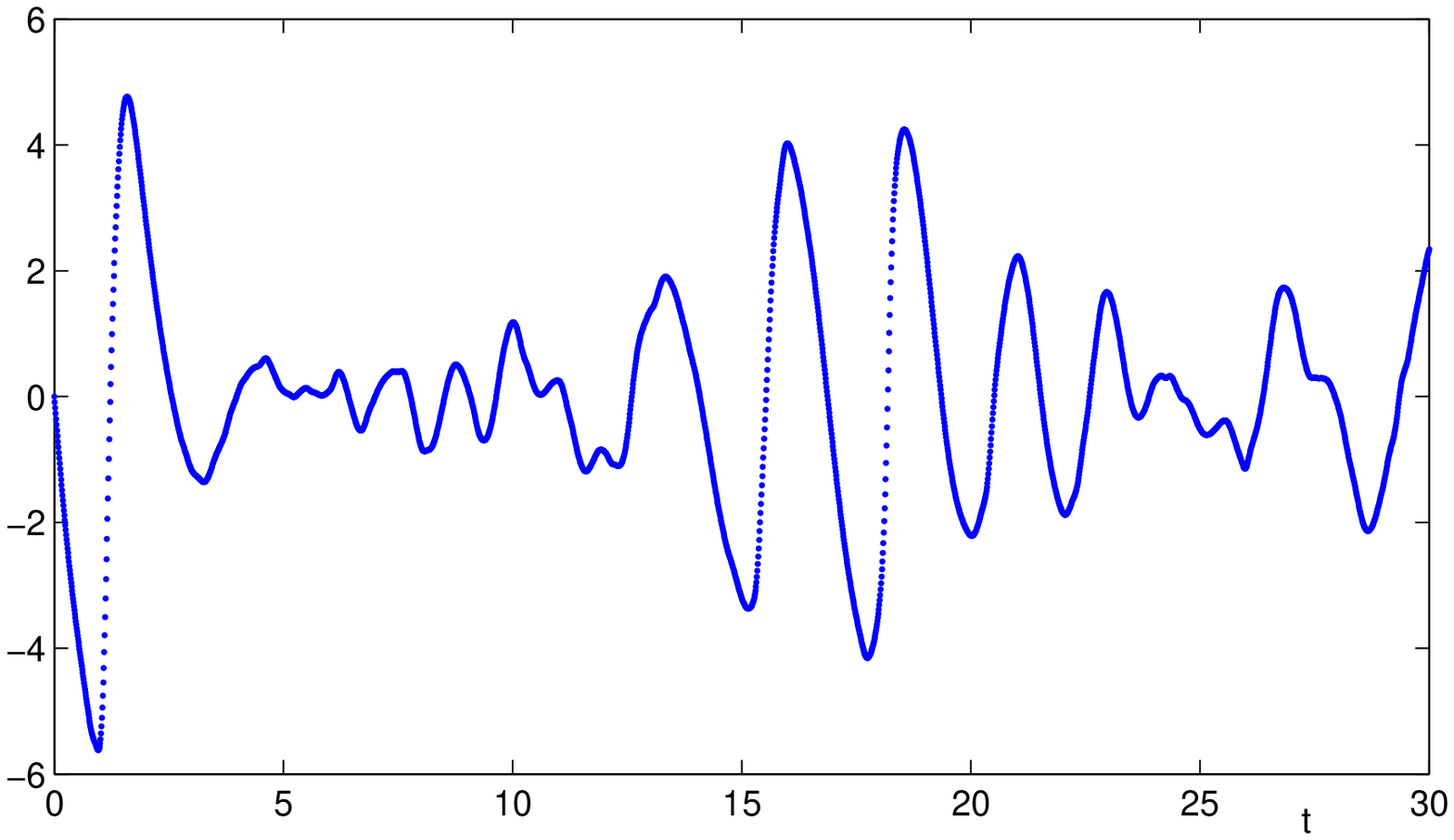}
\includegraphics[scale=0.35]{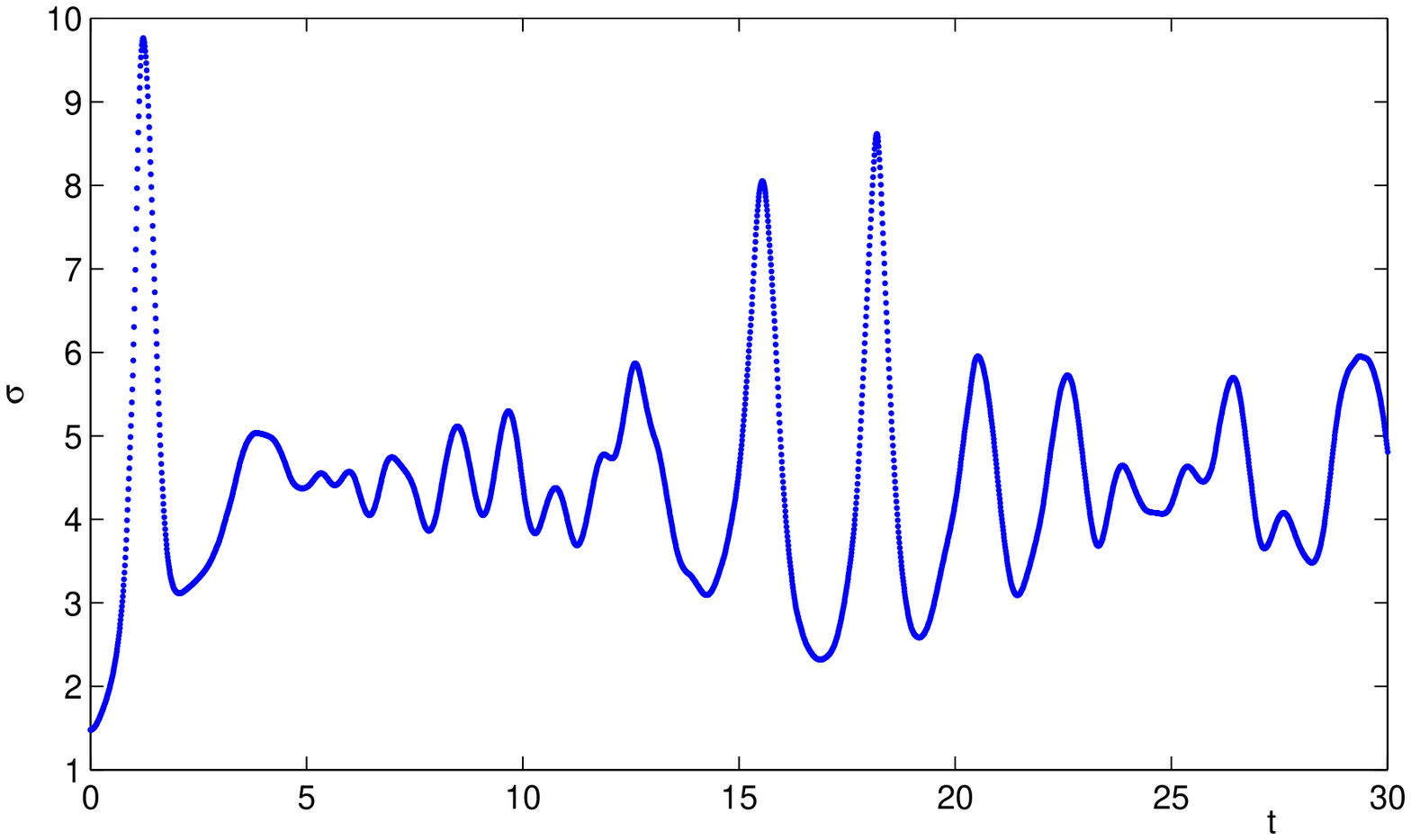}
\includegraphics[scale=0.4]{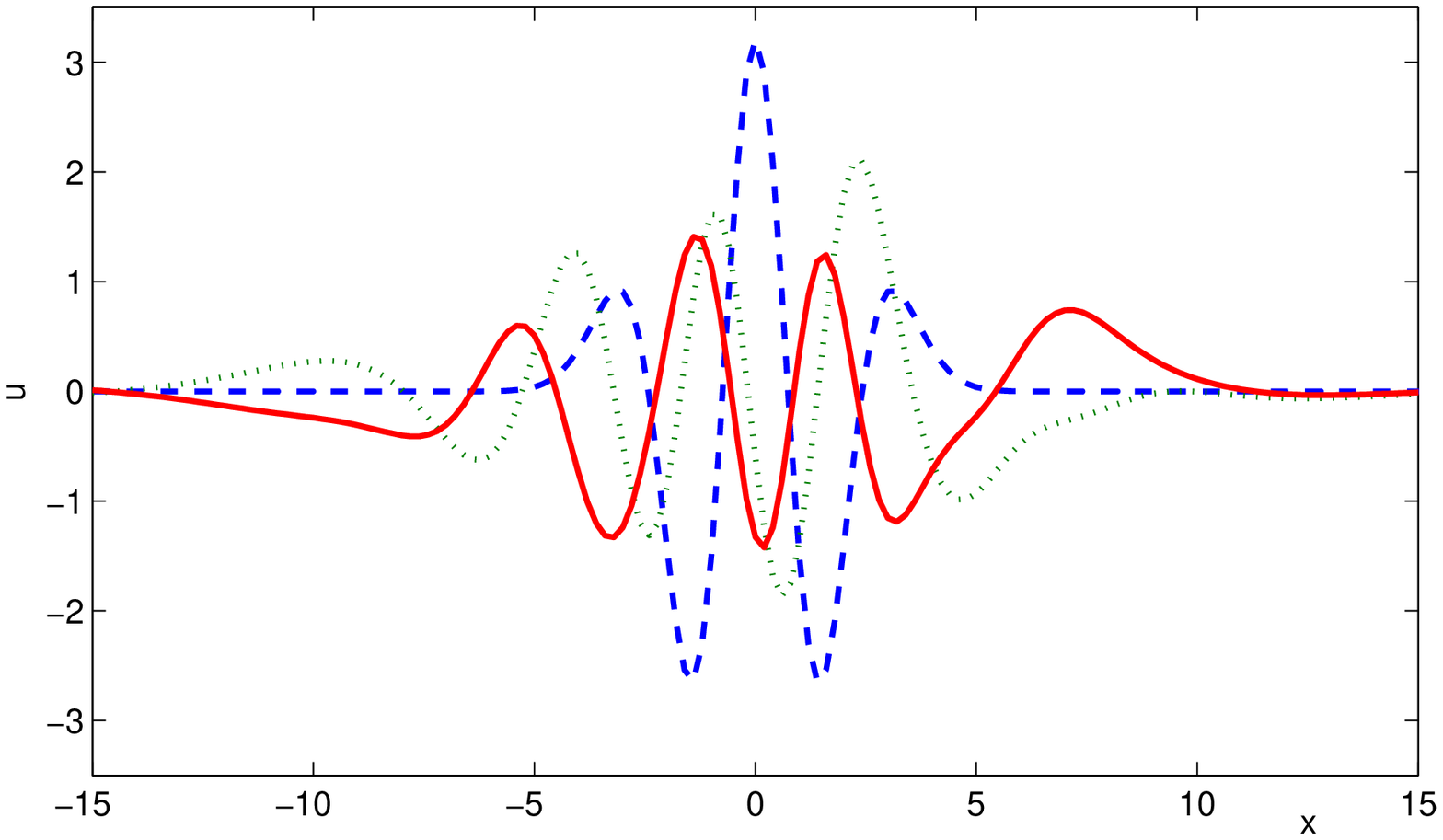}
\end{center}
\caption{\label{fig-Cauchy-1}
Numerical solution of the linearized log--KdV equation
(\ref{linlogKdV}) with the even initial data (\ref{even-data}):
center of mass $\bar{x}$ (top left) and standard deviation $\sigma$ (top right) versus time;
the profile $u(x,t)$ versus $x$ (bottom) for
$t = 0$ (blue dashed), $t = 2.5$ (green dotted), and $t = 5$ (red solid).}
\end{figure}

Figure \ref{fig-Cauchy-1} reports numerical computations for the even initial data
\begin{equation}
\label{even-data}
u_0(x) = \left(x^4 - \frac{3 (3 + 4\alpha)}{(1+\alpha)(1+2\alpha)} x^2
  + \frac{6}{(1+\alpha)(1+2\alpha)} \right) e^{-\frac{(1+2\alpha) x^2}{4}},
\end{equation}
with $\alpha = 0.25$. The even function (\ref{even-data}) also satisfies the symplectic orthogonality conditions
$\langle v_G, u_0 \rangle_{L^2} = \langle \partial_x^{-1} v_G, u_0 \rangle_{L^2} = 0$ for any $\alpha \geq 0$.
The results are qualitatively similar to the case of odd initial data. The preserved localization
of the solution coexists with the radiation at the left slope of the Gaussian pulse.

\subsection{Discussion}

The numerical results of Section 4.2 support the conjecture that if
the initial data $u_0 \in X_c$ of the linearized log--KdV equation
(\ref{linlogKdV}) satisfies the constraint
$$
u_0(x) = v_G(x) h_0(x) \quad \mbox{\rm with} \quad
h_0 \in L^{\infty}(\mathbb{R}) \cap L^2(\mathbb{R}),
$$
then there exists a unique solution of the linearized log--KdV equation
(\ref{linlogKdV}) in the form
$$
u(x,t) = v_G(x) h(x,t),
$$
where $h(\cdot,t) \in L^{\infty}(\mathbb{R}) \cap L^2(\mathbb{R})$ for every $t \in \mathbb{R}$ and
$h(\cdot,0) = h_0$.

If this result can be established rigorously, then one can analyze the
linearized log--KdV equation with a source term to complete
analysis of the equivalent log--KdV equation (\ref{logKdV-w}) with the initial
data $w_0 = v_G h_0$, where $h_0 \in L^{\infty}(\mathbb{R}) \cap L^2(\mathbb{R})$
and $\| h_0 \|_{L^{\infty} \cap L^2}$ is sufficiently small.
This route may lead to the proof of
orbital stability of the Gaussian solitary wave $v_G$ for the log--KdV equation (\ref{logKdV}) in the class
of functions with initial data $v_0 = (1 + h_0) v_G$, where $h_0 \in L^{\infty}(\mathbb{R}) \cap L^2(\mathbb{R})$
and $\| h_0 \|_{L^{\infty} \cap L^2}$ is sufficiently small.
However, this work is still to be done, hence the problem of nonlinear orbital stability of the
Gaussian solitary wave $v_G$ remains opened for further studies.

\vspace{0.5cm}

{\bf Acknowledgement.} D.P. is indebted to Guillaume James (University of Grenoble)
for bringing up the problem and fruitful collaboration, as well as to Thierry Gallay
(University of Grenoble) who noticed the visible radiation in our numerical
simulations. D.P. is supported by the CNRS Visiting Fellowship. He thanks
Institut de Math\'ematiques et de Mod\'elisation,
Universit\'e Montpellier for hospitality and CNRS for support during his visit
(September-November, 2013).

\bibliographystyle{amsplain}
\bibliography{biblio}

\providecommand{\bysame}{\leavevmode\hbox to3em{\hrulefill}\thinspace}
\providecommand{\MR}{\relax\ifhmode\unskip\space\fi MR }
\providecommand{\MRhref}[2]{%
  \href{http://www.ams.org/mathscinet-getitem?mr=#1}{#2}
}
\providecommand{\href}[2]{#2}
\begin{thebibliography}{10}

\bibitem{Pava}
J.~Angulo~Pava, \emph{Nonlinear dispersive equations}, Mathematical Surveys and
  Monographs, vol. 156, American Mathematical Society, Providence, RI, 2009,
  Existence and stability of solitary and periodic travelling wave solutions.

\bibitem{CG}
R.~Carles and C.~Gallo, \emph{Finite time extinction by nonlinear damping for
  the {S}chr{\"o}dinger equation}, Comm. Part. Diff. Eq. \textbf{36} (2011),
  no.~6, 961--975.

\bibitem{CaOz-p}
R.~Carles and T.~Ozawa, \emph{Finite time extinction for nonlinear
  {S}chr{\"o}dinger equation in {1D} and {2D}}, preprint, archived at
  \url{http://arxiv.org/abs/1405.0995}, 2014.

\bibitem{Caz2}
T.~Cazenave, \emph{Stable solutions of the logarithmic {S}chr\"odinger
  equation}, Nonlinear Anal. \textbf{7} (1983), no.~10, 1127--1140.

\bibitem{Caz}
T.~Cazenave, \emph{Semilinear {S}chr\"odinger equations}, Courant Lecture Notes
  in Mathematics, vol.~10, New York University Courant Institute of
  Mathematical Sciences, New York, 2003.

\bibitem{Caz1}
T.~Cazenave and A.~Haraux, \emph{\'{E}quations d'\'evolution avec non
  lin\'earit\'e logarithmique}, Ann. Fac. Sci. Toulouse Math. (5) \textbf{2}
  (1980), no.~1, 21--51.

\bibitem{Chat}
A.~Chatterjee, \emph{Asymptotic solution for solitary waves in a chain of
  elastic spheres}, Phys. Rev. E \textbf{59} (1999), 5912--5919.

\bibitem{DS}
N.~Dunford and J.~T. Schwartz, \emph{Linear operators. {P}art {II}: {S}pectral
  theory. {S}elf adjoint operators in {H}ilbert space}, With the assistance of
  William G. Bade and Robert G. Bartle, Interscience Publishers John Wiley \&
  Sons\ New York-London, 1963.

\bibitem{GV}
J.~Ginibre and G.~Velo, \emph{The global {C}auchy problem for the nonlinear
  {S}chr\"odinger equation revisited}, Ann. Inst. H. Poincar\'e Anal. Non
  Lin\'eaire \textbf{2} (1985), 309--327.

\bibitem{JP13}
G.~James and D.~Pelinovsky, \emph{Gaussian solitary waves and compactons in
  {F}ermi-{P}asta-{U}lam lattices with {H}ertzian potentials}, Proc. Roy. Soc.
  A \textbf{470} (2014), 20130465 (20 pages).

\bibitem{kapstef}
T.~Kapitula and A.~Stefanov, \emph{A {H}amiltonian--{K}rein (instability) index
  theory for {KdV}-like eigenvalue problems}, Stud. Appl. Math. \textbf{132}
  (2014), 183--211.

\bibitem{Kato}
T.~Kato, \emph{On the {K}orteweg-de\thinspace {V}ries equation}, Manuscripta
  Math. \textbf{28} (1979), no.~1-3, 89--99.

\bibitem{KPV93}
C.~E. Kenig, G.~Ponce, and L.~Vega, \emph{Well-posedness and scattering results
  for the generalized {K}orteweg-de {V}ries equation via the contraction
  principle}, Comm. Pure Appl. Math. \textbf{46} (1993), no.~4, 527--620.

\bibitem{Miller}
P.~D. Miller, \emph{Applied asymptotic analysis}, Graduate Studies in
  Mathematics, vol.~75, American Mathematical Society, Providence, RI, 2006.

\bibitem{Nesterenko}
V.~F. Nesterenko, \emph{Dynamics of heterogeneous materials}, Springer Verlag,
  New York, 2001.

\bibitem{dmitrystab}
D.~E. Pelinovsky, \emph{Spectral stability of nonlinear waves in {KdV}-type
  evolution equations}, Spectral analysis, stability, and bifurcation in modern
  nonlinear physical systems (O.~Kirillov and D.~E. Pelinovsky, eds.),
  Wiley--ISTE, 2014, pp.~377--400.

\bibitem{Teschl}
G.~Teschl, \emph{Ordinary differential equations and dynamical systems},
  Graduate Studies in Mathematics, vol. 140, American Mathematical Society,
  Providence, RI, 2012.

\bibitem{Zhidkov}
P.~E. Zhidkov, \emph{Korteweg-de {V}ries and nonlinear {S}chr\"odinger
  equations: qualitative theory}, Lecture Notes in Mathematics, vol. 1756,
  Springer-Verlag, Berlin, 2001.

\end{thebibliography}
\end{document}